\documentclass{amsart}

\usepackage{amssymb,latexsym,amscd}

\numberwithin{equation}{section}

\newtheorem{theorem}{Theorem}[section]
\newtheorem{proposition}[theorem]{Proposition}
\newtheorem{corollary}[theorem]{Corollary}

\newtheorem{lemma}[theorem]{Lemma}

\def\ii{\mathbf{i}}
\def\AA{\mathbb{A}}

\def\GG{\mathbb{G}}
\def\QQ{\mathbb{Q}}
\def\RR{\mathbb{R}}
\def\ZZ{\mathbb{Z}}

\def\gg{\mathfrak{g}}

\begin{document}

\title{Geometric and unipotent crystals}

\author{Arkady  Berenstein}

\author{David Kazhdan}

\address{Department of Mathematics, Harvard University}
\email{arkadiy@math.harvard.edu}


\address{Department of Mathematics, Harvard University} 
\email{kazhdan@math.harvard.edu}

\date{December 10, 1999}

\thanks{The research of both authors was supported in part
by NSF grants}

\maketitle

\makeatletter
\renewcommand{\@evenhead}{\tiny \thepage \hfill  
A.~BERENSTEIN and D.~KAZHDAN \hfill}

\renewcommand{\@oddhead}{\tiny \hfill  
GEOMETRIC AND UNIPOTENT CRYSTALS
 \hfill \thepage}
\makeatother

\tableofcontents


\section{Introduction}

\label{sect:introduction} 

Let $G$ be a split semisimple algebraic group over $\QQ$, 
$\gg$ be the Lie algebra of $G$ and  $U_q(\gg)$ be the corresponding  
{\it quantized enveloping algebra}. Lusztig has introduced in \cite{lu90} 
{\it canonical bases}  for finite-dimensional $U_q(\gg)$-modules. About the same time  
Kashiwara introduced in \cite{k90}  {\it crystal bases} 
as a natural framework for parametrizing bases of finite-dimensional 
$U_q(\gg)$-modules. It was shown in \cite{lu92} that Kashiwara's 
crystal bases are the limits as $q\to 0$ of Lusztig's canonical bases.  
Later, in \cite{k93} Kashiwara introduced a new combinatorial concept: 
{\it crystals}. Kashiwara's crystals generalize the crystal bases and 
provide a natural framework for their study.

In this paper we study  {\it geometric crystals} 
and {\it unipotent crystals} which are algebro-geometric analogues of
the crystals and the crystal bases respectively. 

Our approach of the ``geometrization'' of combinatorial 
objects comes back to the works of Lusztig. On the one hand, Lusztig constructed  
for each reduced  decomposition $\ii$ of the longest 
element $w_0$ of the Weyl group $W$
of $G$ a parametrization $\psi^{\ii}$ of the canonical basis ${\bf
B}$ by the cone  
$(\ZZ_{\ge 0})^{l_0}$, where $l_0$ is the length $w_0$. 
On the other hand,  in his study of the positive elements of $G(\RR)$ 
Lusztig constructed  for each  $\ii$ a  birational isomorphism 
$\pi^{\ii}:\AA^{l_0}\widetilde \to  U$, where $U$ is the unipotent radical of the 
Borel subgroup of $G$.
Moreover, Lusztig observed that for any two decompositions 
$\ii,\ii'$ of $w_0$,
the piecewise-linear transformation $(\psi^{\ii'})^{-1}\circ \psi^{\ii}$
of $(\ZZ_{\ge 0})^{l_0}$ is similar to  the  birational transformation
$(\pi^{\ii'})^{-1}\circ \pi^\ii$ of $\AA^{l_0}$
(see \cite{lu},\,\cite{bfz}). 

This analogy was further developed in a series of papers 
\cite{bfz},\cite{bz2},\,\cite{bz3},\,\cite{fz}. 
The main idea of these papers is to introduce the notion of the
``tropicalization" -- a set of  formal rules for
the passage from birational isomorphisms to piecewise-linear maps, 
and to study piecewise-linear automorphisms of ${\bf B}$ using the
``tropicalization" of explicit birational isomorphisms
$\pi^\ii:\AA^{l_0}\widetilde\to U$.  

Another geometric approach to constructing crystal bases is 
suggested in the recent work of Braverman and Gaitsgori \cite{bg}. Within this approach 
the crystal 
bases and, in fact,  actual bases in finite-dimensional $G$-modules,
are constructed in terms of perverse sheaves on the 
affine Grassmannian of $G$.

In the present paper we go in the opposite direction -- we study 
geometric crystals in their own right and ``geometrize" some of the 
mentioned above piecewise-linear structures of the crystal bases. 
 
Let $I$ be the set of vertices of the Dynkin diagram of $G$, and $T$ be
a maximal torus of $G$. The structure of a 
 {\it geometric crystal} on an algebraic variety $X$
consists of a rational morphism $\gamma :X\to T$ and a compatible family 
$e_i:\GG_m\times X\to X$, $i\in I$,  of  rational actions of the
multiplicative group  $\GG_m$ on $X$. We show that such a structure
induces a rational action of the Weyl group $W$ on $X$. Surprisingly, 
many interesting rational actions of $W$ come from 
geometric crystals. For example, the natural action of 
$W$ on Grothendieck's simultaneous resolution 
$\tilde G \to G$ comes from a structure of a geometric crystal 
on $\tilde G$. Also all the examples of the action of  
$W$  in \cite {bk}, come from geometric crystals. Another application of
geometric crystals is a construction for any  $SL_n$-crystal built on 
$(X,\gamma)$ of a  {\it trivialization}: a 
$W$-equivariant isomorphism $X\widetilde \to \gamma ^{-1}(e) \times T$.

It is also interesting that the Langlands dual group  
${}^LG$ emerges when we reconstruct  
Kashiwara's combinatorial crystals out of {\it positive} 
geometric crystals (see section \ref{subsect:tropicalization}). 
The presence of ${}^LG$ in the ``crystal world" 
has been noticed in \cite{lu92}. The combinatorial 
results of \cite{bz3} also involve ${}^LG$.

A number of statements in the paper are almost immediate. We call them
{\it Lemmas}. Proofs of all the Lemmas from sections
\ref{sect:definitions},  \ref {sect:unipotent crystals} 
and \ref{sect:examples} are left for the reader. 

\smallskip

The material of the paper is organized as follows.

\smallskip
$\bullet$ In section \ref{sect:definitions} 
we introduce geometric crystals and formulate their main properties. 

\smallskip
$\bullet$ In section \ref {sect:unipotent crystals} we introduce 
{\it unipotent crystals} as an algebro-geometric analogue of 
crystal bases for $U_q({}^L\gg)$-modules, where  ${}^L\gg$ is the Lie 
algebra of ${}^L G$.  A unipotent crystal on an 
irreducible algebraic variety $X$ consists
of a rational action of $U$ on $X$ and a rational $U$-equivariant morphism 
${\bf f}:X\to G/U$.
The unipotent crystals have a number of properties
characteristic of the crystal bases.
One of the main properties is a natural {\it product} of unipotent crystals 
which is an analogue of tensor product of crystal bases. We prove 
that the category of unipotent crystals is strict monoidal with respect to 
this product. We also define the notion of a {\it dual unipotent crystal}. 
One of the main results of this section is a
construction for any unipotent crystal on $X$, of an {\it induced}
geometric crystal ${\mathcal X}$. We also give a closed formula 
for the action of $W$ on this induced geometric crystal 
${\mathcal X}$. 

\smallskip
$\bullet$ In section \ref{sect:examples} we study the 
{\it standard} unipotent crystals which are the 
$U$-orbits $BwB/B\subset G/B$. We obtain a decomposition 
of any standard crystal  $BwB/B$ into the product of 
$1$-dimensional crystals corresponding to any reduced decomposition 
of $w\in W$.  This is a geometrization of Kashiwara's results. 
We also construct $W$-invariant functions 
on certain standard unipotent crystals. These functions are important for
the study of $\gamma$-functions of representations ${}^LG$ (see \cite{bk}).
\smallskip

$\bullet$ In section \ref{sect:proofs} we collect proofs of the 
results from sections \ref{sect:definitions},  
\ref {sect:unipotent crystals} and \ref{sect:examples}. 
\smallskip

$\bullet$ In section \ref{sect:projections} we study the 
restrictions to a standard Levi subgroup 
$L\subset G$ of standard unipotent $G$-crystals $BwB/B$. In particular, 
for each $w\in W$ we construct a rational morphism 
${\bf p}_w:BwB/B\to L/(U\cap L)$. 
In the case when 
${\bf p}_w$ is a birational isomorphism with its image, we 
study the direct image of $W_L$-invariant functions under ${\bf p}_w$. 
In the special case when $G=GL(m+n),L=GL(m)\times GL(n)$, and  
$w=w_0^L\cdot w_0$, it 
is possible to  write explicitly 
the corresponding action of the group $W_L$  on the image of 
${\bf p}_w$, and the $W_L$-equivariant trivialization. This 
simplification of the 
structure of the crystal on $BwB/B$ is used in \cite {bk} for the 
prove that Piatetsky-Shapiro's $\gamma$-function to
$GL(m)\times GL(n)$ is equal to the one introduced in \cite {bk}. 
We expect that in general the action of $W_L$ on $BwB/B$, and the 
trivialization of the corresponding crystal would be easy to describe, 
and that this description can help in the study of 
$\gamma$-functions on $L$.

\smallskip

$\bullet$ In the Appendix we collect necessary definitions and results related
to Kashiwara's crystals. We refer to these crystals as 
{\it combinatorial} in order to distinguish between 
them and their geometric counterparts.  

\smallskip
The notion of  {\it a geometric crystal} was introduced by the  first 
author who noticed that the formulas for the $W$-action which appear in 
the definition of $\gamma$-functions for the group $GL_2 \times GL_3$ 
(see \cite {bk}) can be interpreted as a ``geometrization" of the 
$W$-action for the free combinatorial crystal defined earlier by the 
first author. 
We want to express our gratitude to Alexander Braverman for his help in different 
stages of this work and in particular for the definition of the notion 
of {\it degree} in section \ref{subsect:tropicalization}. 
We are also grateful to Yuval Flicker for numerous remarks about the paper. 
Special thanks are due to Charles Cochet for correcting errors in the first version of the paper.

\section {Definitions and  main results on geometric crystals}
\label{sect:definitions}

\subsection {General notation} 

\label{subsect:notation}

Let $G$ be a split reductive algebraic group over $\QQ$ and 
$T\subset G$ a maximal torus.  We denote by  $\Lambda^\vee$ and 
$\Lambda$ the
lattices of  co-characters and characters of $T$ and by  
$ \left<\cdot,\cdot\right>$  
the evaluation pairing $\Lambda^\vee\times \Lambda \to \ZZ$.

Let $B$ be a  Borel subgroup containing $T$. Denote by  $I$  the set of 
vertices of the Dynkin diagram of $G$; for any $i\in I$ denote by 
$\alpha_i\in \Lambda$ the simple root $\alpha_i:T\to \GG_m$, 
and by $\alpha_i^\vee\in \Lambda^\vee$ the simple coroot 
$\alpha_i^\vee:\GG_m\to T$. 

Let $B^-\subset G$ be the Borel subgroup containing $T$ such that  
$B\cap B^-=T$.  Denote by $U$ and $U^-$ 
respectively the unipotent radicals of $B$ and $B^-$.

For each $i\in I$ we denote by 
$U_i\subset U$ and $U^-_i\subset U^-$ the corresponding simple root
subgroups and denote 
by $\xi_i:U\to U_i, \xi^-_i:U^-\to U^-_i$ 
the canonical projections. 

We fix a family of isomorphisms $x_i: \GG_a\widetilde \to U_i,y_i:\GG_a\widetilde\to U^-_i$, 
$i\in I$, such that 
\begin{equation}
\label{eq:commutation}
x_i(a)y_i(a')=y_i\left(\frac{a'}{1+aa'}\right)\alpha_i^\vee(1+aa')
x_i\left(\frac{a}{1+aa'}\right) \ . 
\end{equation}
Clearly, each isomorphism $y_i$ is uniquely determined by  $x_i$. 

\medskip

\noindent {\bf Remark}. Each pair $x_i,y_i$ defines a homomorphism 
$\phi_i:SL_2\to G$:
$$\phi_i\begin{pmatrix} a & b \\
                c & d 
\end{pmatrix}=y_i\left(\frac{c}{a}\right)\alpha_i^\vee(a)
x_i\left(\frac{b}{a}\right)$$

\smallskip

Denote by $\widehat U$ and $\widehat U^-$ 
respectively the spaces $Hom(U,\GG_a)$ and $Hom(U^-,\GG_a)$.  
For each $i\in I$ define $\chi_i\in \widehat U, \chi_i^-\in \widehat U^-$ by  
$$\chi_i(u_+)=x_i^{-1}(\xi_i(u_+)),
~\chi^-_i(u_-)=y_i^{-1}(\xi^-_i(u_-))$$ 
for $u_-\in U^-$, $u_+\in U$. By definition, the family 
$\chi_i$, $i\in I$, is a basis in the vector space $\widehat U$, 
and the family $\chi_i^-$, $i\in I$, is a basis in the vector 
space $\widehat U^-$. For any $\chi\in \widehat U$, and any 
$\chi^-\in \widehat U^-$ 
define functions $\overline \chi,\overline \chi^{\,-}:U^-\cdot T\cdot U
\to \GG_a$ by 
$$\overline \chi^{\,-}(u_-tu_+)=\chi^-(u_-), 
~\overline \chi(u_-tu_+)=\chi(u_+) \ .$$

The Weyl group  $W=Norm_G(T)/T$ of $G$ is generated by 
simple reflections $s_i\in W, i\in I$. The group 
$W$ acts on the lattices $\Lambda$, $\Lambda^\vee$ by
$s_i(\lambda)=\lambda-
\left<\alpha_i^\vee,\lambda\right>\alpha_i $ for 
$\lambda\in \Lambda$,  
$s_i(\lambda^\vee)=\lambda^\vee-
\left<\lambda^\vee,\alpha_i\right>\alpha_i^\vee$ for 
$\lambda^\vee\in \Lambda^\vee$. 

Let $l:W\to \ZZ_{\ge 0}$ ($w\mapsto l(w)$) be the length function. For any  
sequence $\ii=(i_1,\ldots,i_l)\in I^l$ we write
$w(\ii)=s_{i_1}\cdots s_{i_l}$. 
A sequence $\ii\in I^l$ is called {\it reduced}  
if  the length of $w(\ii)$ is equal to $l$. 
For any $w\in W$ we denote by $R(w)$ 
the set of all reduced sequences $\ii$ such that $w(\ii)=w$. We denote
  by  $w_0\in W$ the element of the maximal length in $W$.

For $i\in I$ define 
$\overline {s_i}\in G$ by
\begin{equation}
\label{eq:split si}
\overline s_i=x_i(-1)y_i(1)x_i(-1) =\phi_i
\begin{pmatrix} 
0 & -1 \\
1 & 0 
\end{pmatrix} \ .
\end{equation}

Each $\overline s_i$ belongs to $Norm_G(T)$ and is a representative of
$s_i\in W$. It is  well-known (\cite {bou}) that the elements 
$\overline {s_i}$, $i\in I$, satisfy the braid relations. 
Therefore we can associate to each $w\in W$ its 
{\it standard representative} $\overline  w\in Norm_G(T)$ in such a way
that for any 
$(i_1,\ldots,i_l)\in R(w)$ we have: 
\begin{equation}
\label{eq:split w}
\overline w=\overline {s_{i_1}}\cdots \overline {s_{i_l}} \ .
\end{equation}

\subsection{Geometric pre-crystals and geometric crystals}

\label{subsect:geometric crystals}

Let $X$ and $Y$ be algebraic varieties over $\QQ$. 
Denote by ${\bf R}(X,Y)$ the set of all rational morphisms from $X$ to $Y$. 

For any ${\bf f}\in {\bf R}(X,Y)$ denote by  $dom({\bf f})\subset X$ the 
maximal open subset of $X$ on which $f$ is defined; denote by 
${\bf f}_{reg}:dom({\bf f})\to Y$ the corresponding regular morphism. We
denote by  $ran({\bf f})\subset Y$ the closure of the constructible set
${\bf f}_{reg}(dom({\bf f}))$ in $Y$. 
For any regular morphism $f:X'\to Y$, where $X'\subset X$ is a 
dense subset we denote by $[f]:X\to Y$ 
the corresponding rational morphism. Note that $[{\bf f}_{reg}]={\bf f}$ and 
$dom({\bf f}_{reg})=dom({\bf f})$ for any ${\bf f}\in {\bf R}(X,Y)$.

It is easy to see that for any
irreducible  algebraic varieties $X,Y,Z$ and
rational morphisms  $f:X\to Y,~g:Y\to Z$ such that 
$dom(g)$ intersects $ran(f)$ non-trivially, the composition $(f,g)\mapsto 
g\circ f$ is well-defined and is a  rational morphism from $X$ to $Z$.

We denote by  ${\mathcal V}$ be the category whose 
objects are irreducible algebraic varieties and 
arrows are dominant rational morphisms.

For any algebraic group $H$ we call a rational action 
$\alpha:H\times X\to X$ {\it unital} if 
$dom(\alpha)\supset \{e\}\times X$.

\noindent {\bf Definition}. Let $X\in Ob({\mathcal V})$ and 
$\gamma$ be a rational morphism $X\to T$.  
A {\it geometric $G$-pre-crystal} (or simply a {\it geometric pre-crystal}) 
on $(X,\gamma)$, is a family 
$e_i:\GG_m\times X\to X$, $i\in I$,  of unital rational 
actions of the multiplicative group 
$\GG_m$: $(c,x)\mapsto e_i^c(x)$, such that 
$\gamma(e_i^c(x))=\alpha_i^\vee(c)\gamma(x)$.

\medskip

\noindent {\bf Remark}. Geometric pre-crystals are analogues of   
free combinatorial pre-crystals (see Appendix). 
Under this analogy, the 
variety $X$ corresponds to the set $B$,  
the maximal torus $T$ corresponds to the lattice $\Lambda^\vee$, 
the rational morphism $\gamma:X\to T$ corresponds to 
the map $\tilde \gamma:B\to \Lambda^\vee$ and rational actions $e_i$ of  
$\GG_m$ on $X$ corresponds to bijections  $\tilde e_i:B\to B$. 
We will make this analogy precise in section 
\ref{subsect:tropicalization}.

\smallskip

Given a geometric pre-crystal ${\mathcal X}$, and a reduced sequence 
$\ii=(i_1,\ldots,i_l)\in I^l$, we define a rational morphism 
$e_{\ii}:T\times X \to X$  by 
\begin{equation}
\label{eq:general eii}
(t,x)\mapsto e_\ii^t(x)=e_{i_1}^{\alpha^{(1)}(t)}\circ 
\cdots \circ e_{i_l}^{\alpha^{(l)}(t)}(x) \ ,
\end{equation}
where $\alpha^{(k)}=s_{i_l}s_{i_{l-1}}\cdots s_{i_{k+1}}(\alpha_{i_k})$, $k=1,\ldots,l$ 
are the {\it associated} positive roots. 

\noindent {\bf Definition}. 
A geometric pre-crystal ${\mathcal X}$ on $(X,\gamma)$ is called a
{\it geometric crystal} 
if for any $w\in W$,  and any $\ii,\ii'\in R(w)$,  one has:
\begin{equation}
\label{eq:verma}
e_\ii=e_{\ii'} \ .
\end{equation}

\begin{lemma}  
\label{le:verma}
The relations (\ref{eq:verma}) are equivalent to the following relations 
between $e_i, e_j$ for $i,j\in I$:
\begin{equation}
\label{eq:A1A1}\displaystyle{e_i^{c_1} e_j^{c_2}=e_j^{c_2} e_i^{c_1}}
\end{equation}
if $\left<\alpha_i^\vee,\alpha_j\right>=0$;  
\begin{equation}
\label{eq:A2}\displaystyle{e_i^{c_1}e_j^{c_1c_2}e_i^{c_2}
=e_j^{c_2}e_i^{c_1c_2}e_j^{c_1}}
\end{equation}
if $\left<\alpha_i^\vee,\alpha_j\right>
=\left<\alpha_j^\vee,\alpha_i\right>=-1$;
\begin{equation}
\label{eq:B2}\displaystyle{e_i^{c_1}e_j^{c_1^2c_2}e_i^{c_1c_2}e_j^{c_2}
=e_j^{c_2}e_i^{c_1c_2}e_j^{c_1^2c_2}e_i^{c_1}}
\end{equation}
if  $\left<\alpha_i^\vee,\alpha_j\right>=-2, 
\left<\alpha_j^\vee,\alpha_i\right>=-1$;
\begin{equation}
\label{eq:G2}\displaystyle{
e_i^{c_1}
e_j^{c_1^3c_2}
e_i^{c_1^2c_2}
e_j^{c_1^3c_2^2}
e_i^{c_1c_2}
e_j^{c_2}
=e_j^{c_2}
e_i^{c_1c_2}
e_j^{c_1^3c_2^2}
e_i^{c_1^2c_2}
e_j^{c_1^3c_2}
e_i^{c_1}}
\end{equation}
if $\left<\alpha_i^\vee,\alpha_j\right>=-3, 
\left<\alpha_j^\vee,\alpha_i\right>=-1$. 

\end{lemma}

\begin{lemma} Let ${\mathcal X}$ be a geometric pre-crystal.
For any $\lambda\in \Lambda$, $i\in I$,
the formula $x\mapsto e_i^{\lambda(\gamma(x))}(x)$ 
defines a rational morphism $X\to X$. This is a birational isomorphism $X\widetilde\to X$ if and only if $\left< \alpha_i^\vee,\lambda\right>\in \{-2,0\}$. 
In the latter case the inverse morphism is given by the formula 
$$x\mapsto 
\begin{cases}
e_i^{\lambda(\gamma(x))}(x),
& \text{if $\left<\alpha_i^\vee,\lambda\right>=-2$,} \\
e_i^{\frac{1}{\lambda(\gamma(x))}}(x), &
\text{if $\left<\alpha_i^\vee,\lambda\right>=0$.}
\end{cases} $$

\end{lemma}

\noindent {\bf Remarks}. 

\noindent 1. The relations (\ref{eq:A1A1})--(\ref{eq:G2}) are 
multiplicative 
analogues of the Verma relations in the universal 
enveloping algebra $U({}^L\gg)$ (see \cite{lu}, Proposition 39.3.7). 

\noindent 2. An analogue of the relations 
(\ref{eq:A1A1}), (\ref{eq:A2}) for a combinatorial $GL_n$-crystal 
was considered in \cite{bkir}, Theorem 1.1.

\smallskip

Let ${\mathcal X}$ be a geometric crystal. For each $w\in W$   
we define a rational morphism 
$e_w:T\times X\to X$  by the formula
$e_w :=e_\ii$ for any  $\ii\in R(w)$ (see (\ref{eq:verma})).

\begin{proposition}
\label{pr:W-action} 
The correspondence $W\times X\to X$ defined by 
\begin{equation}
\label{eq:W-action}
(w,x)\mapsto w(x)=e_w^{\gamma(x)^{-1}}(x)
\end{equation}
is a rational unital action of $W$ on $X$. 
\end{proposition}

\noindent {\bf Remark}. The formula   
$s_i(x)=e_i^{\frac{1}{\alpha_i(\gamma(x))}}(x)$  is 
a multiplicative analogue of (\ref{eq:Kashiwara W}) in the Appendix.

\subsection {The geometric crystal ${\mathcal X}_{w_0}$}

\label{subsect:Xinfty} 

Let $B^-_{w_0}:=U\overline {w_0}U\cap B^-$. The natural inclusion 
$B^-_{w_0}\hookrightarrow G$ induces the open  inclusion ${\bf j}_0:B^-_{w_0}\hookrightarrow 
G/B$. 
Let $\gamma:G/B\to T$ be the rational morphism defined by 
$\gamma=pr_T\circ [{\bf j}_0]^{-1}$, where $pr_T:B^-\to T=B^-/U^-$ is 
the natural projection.

For $i\in I$ let $\varphi_i:G/B\to \GG_a$ be the rational function given 
by  $\varphi_i:=\overline \chi_i^{\,-}\circ [{\bf j}_0]^{-1}$,  
where  $\overline \chi_i^{\, -}:B^-\cdot U\to \GG_a$ 
is the regular function defined in section \ref{subsect:notation}. 

\begin{lemma}
\label{le:properties of phi} For each $i\in I$ we have:

\noindent (a)  $\varphi_i\not \equiv 0$;

\noindent (b) $\frac{1}{\varphi_i(x_i(a)\cdot x)}=
\frac{1}{\varphi_i(x)}+a$;

\noindent (c) $\gamma(x_i(a)\cdot x)=
\alpha_i^\vee(1+a\varphi_i(x))\gamma(x)$.

\end{lemma}

For each $i\in I$ define a rational morphism  
$e_i:\GG_m\times G/B\to G/B$ by the formula  
\begin{equation}
\label{eq:action of ei}
e_i^c(x)= x_i\left(\frac{c-1}{\varphi_i(x)}\right)\cdot x \ .
\end{equation}

\noindent {\bf Remark}. Lemma \ref{le:properties of phi} 
implies the equality
$\varphi_i(e_i^c(x))=c^{-1}\varphi_i(x)$ for $i\in I$.

\begin{theorem} 
\label{th:X(infinity)}
The morphisms $e_i$, $i\in I$,  define a 
geometric crystal on $(G/B,\gamma)$
\end{theorem}

We denote this crystal by ${\mathcal X}_{w_0}$. 
This is one of our main examples of geometric crystals.

\smallskip

\noindent {\bf Remark}. The geometric crystal ${\mathcal X}_{w_0}$ is an 
analogue of the free combinatorial $W$-crystal ${\mathcal B}_{w_0}$,   
and the rational functions $\varphi_i:G/B\to \GG_m$ are 
analogues of the functions $-\tilde \varphi_i$ on ${\mathcal B}_{w_0}$ 
(see Appendix).

\subsection {Positive geometric crystals and their tropicalization}

\label{subsect:tropicalization}

Let $T'$ be an algebraic torus over $\QQ$. We denote by  
$X^\star(T')$ and  $X_\star(T')$ respectively the the lattices of 
characters and co-characters of $T'$, and by $\left<\cdot,\cdot\right>$ 
the canonical pairing $X_\star(T')\times X^\star(T')\to \ZZ$.
 
Let ${\mathcal L}(T')$ be the set of {\it formal 
loops} $\phi:\GG_m\to T'$. 
In particular, ${\mathcal L}(\GG_m)={\mathcal L}(c)$, where ${\mathcal L}(c)$ is 
the set of all Laurent series in the variable $c$.  

Denote by  ${\mathcal L}_0(T')$  the set of {\it formal
disks}, that is,  the set of all $\phi\in {\mathcal L}(T')$ such that 
for any $\mu\in X^\star(T')$ we have $\mu\circ \phi\in {\mathcal
L}_0(c)$, where 
${\mathcal L}_0(c)\subset {\mathcal L}(c)$ is the set of all invertible Taylor 
series in the variable $c$.

Clearly, ${\mathcal L}_0(T')$ has a natural structure of an irreducible  
pro-algebraic variety.

\begin{lemma} The multiplication map 
$X_\star(T')\times {\mathcal L}_0(T')\to {\mathcal L}(T')$
is a bijection. 

\end{lemma}

This defines a surjective map
$deg_{T'}:{\mathcal L}(T')\to X_\star(T')$.

\noindent {\bf Remarks}. 

\noindent 1. If $T'=\GG_m$ then 
 $deg_{\GG_m}:{\mathcal L}(c)\to \ZZ$ is the valuation map which associates 
to a non-zero Laurent series $\phi(c)$ 
its lowest degree. By definition, ${\mathcal L}_0(c)=deg_{\GG_m}^{-1}(0)$.

\noindent 2. For any co-character $\lambda \in X_\star (T')$, the pre-image 
$deg_{T'}^{-1}(\lambda)\subset {\mathcal L}(T')$ is 
naturally isomorphic to ${\mathcal L}_0(T')$. Hence it 
makes sense to talk  about generic points of $deg_{T'}^{-1}(\lambda)$. 

\smallskip

For any $f\in {\bf R}(T',T'')$ and any $\lambda\in X_\star(T'),\mu\in X_\star(T'')$ 
let $U_f(\lambda,\mu)$ be the set of all $\phi\in deg_{T'}^{-1}(\lambda)$ such that 
$f\circ \phi\in deg_{T''}^{-1}(\mu)$.

\begin{lemma} 
\label{le:generic degree} Let $f\in {\bf R}(T',T'')$. For any $\lambda \in X_\star(T')$ 
there is a unique $\mu\in X_\star(T'')$ such that the set $U_f(\lambda,\mu)$ is dense in 
$deg_{T'}^{-1}(\lambda)$. 
\end{lemma}

This allows to define a map $deg(f):X_\star(T')\to X_\star (T'')$ by $deg(f)(\lambda)=\mu$, 
where $\mu\in X_\star(T'')$ is determined in Lemma \ref{le:generic degree}. 

We call this map the {\it degree} of $f$. 
It is easy to check that $deg(f)$ is a piecewise-linear map. 
Note that $deg(f)$ is linear if $f$ is a group homomorphism.

\smallskip 
\noindent {\bf Definition}. 

\noindent (a) A rational function $f$ on a torus $T'$ is
called {\it positive} if it can be written as a ratio $f=f'/f''$, where $f'$
and $f''$ are linear combinations of characters with positive integer coefficients.

\noindent (b) For any two algebraic tori 
$T',T''$ we call a rational morphism $f\in {\bf R}(T',T'')$ 
{\it positive} if for any character 
$\mu:T''\to \GG_m$ the composition 
$\mu\circ f$ 
is a positive rational function on $T'$.

Denote by ${\bf R}_+(T',T'')$ the set of positive rational morphisms from 
$T'$ to $T''$. 

\smallskip

\noindent {\bf Remark}. Any homomorphism of algebraic tori is positive. 

\smallskip

\begin{lemma} 
\label{le:positive} For any two  positive morphisms $f\in {\bf
R}_+(T',T''),g\in {\bf R}_+(T'',T''')$ the composition $g\circ f \in
{\bf R}_+(T',T''')$ is well-defined.
\end{lemma}

Therefore, we can consider a category 
${\mathcal T}_+$ whose objects 
are algebraic tori, and the arrows are 
positive rational morphisms. 

For any $f\in {\bf R}(T',T'')$ and any $\lambda\in X_\star(T'),\mu\in X_\star(T'')$ 
let $U_{f;\lambda}$ be the set of all $\phi\in deg_{T'}^{-1}(\lambda)$ such that 
$f\circ \phi\in {\mathcal L}(T'')$. 

We say that a formal loop $\phi\in {\mathcal L}(T')$ is {\it positive rational} if 
$\phi\in {\bf R}_+(\GG_m,T')$. 

\begin{lemma} 
\label{le:posit} Let $f:T'\to T''$ be a positive morphism. Then 
for any $\lambda \in X_\star (T')$ and any formal positive rational loop 
$\phi\in deg_{T'}^{-1}(\lambda)$ we have
$$deg_{T''}\circ f\circ \phi=deg(f)(\lambda) \ .$$
\end{lemma}

\begin{corollary}
\label{co:composition of degrees}
 For any algebraic tori $T',T'',T'''$ and any 
$f\in {\bf R}_+(T',T'')$, $g\in {\bf R}_+(T'',T''')$ we have 
\begin{equation}
\label{eq:composition of degrees}
deg(g\circ f)=deg(g)\circ deg(f) \ . 
\end{equation} 
\end{corollary}

This implies that there is a functor $Trop:{\mathcal T}_+\longrightarrow 
{\bf Set}$ such that $Trop(T')= X_\star(T')$ and 
$Trop\left(f:T'\to T''\right)=
\left(deg(f):X_\star(T')\to X_\star(T'')\right)$.

Following \cite{bfz}, we call the functor $Trop$ {\it tropicalization}.

\smallskip

\noindent{\bf Remark}.  
If $f:T'\to T''$ is not positive then the assertion of Lemma \ref{le:posit} is not true 
even if $f$ is a birational isomorphism. 
Moreover, one can find a birational isomorphism $f:T'\widetilde\to T''$ such that
(\ref{eq:composition of degrees}) does not hold for the pair $(f,g)$, where $g=f^{-1}:T''\to T'$. 
Consider, for example, the case when $T'=T''=\GG_m,f(c)=c-1,g(c)=c+1$.

\smallskip

\noindent {\bf Definition}. Let ${\mathcal X}$ be a geometric pre-crystal 
on $(X,\gamma)$. A birational isomorphism 
$\theta: T'\widetilde \to X$ for some algebraic 
torus $T'$  is called a {\it positive structure} on
${\mathcal X}$ if the following conditions are satisfied. 

\noindent 1. The rational morphism $\gamma\circ \theta: T'\to T$ is 
positive. 
 
\noindent 2. For each $i\in I$ the rational morphism 
$e_{i;\theta}:\GG_m\times T'\to T'$ 
given by  
\begin{equation}
\label{eq:positive ei}
e_{i;\theta}(c,t')=\theta^{-1}(e_i^{c}(\theta(t')))
\end{equation} 
is positive.

\smallskip

\noindent {\bf Remark}. In all our examples of positive structures $\theta$ 
the composition $\gamma\circ \theta:T'\to T$ is a group homomorphism.

\smallskip

We say that two positive structures $\theta,\theta'$ are {\it equivalent} if
the rational morphisms $\theta^{-1} \circ \theta'$ and  
${\theta'}^{-1} \circ \theta$ are positive.

Recall from the Appendix
that a combinatorial pre-crystal consists of a set $B$, a map 
$\tilde \gamma:B\to \Lambda^\vee$ and  compatible collection of partial 
bijections $\tilde e_i:B\to B$.  

For any  positive structure $\theta$ on a geometric pre-crystal 
${\mathcal X}$ built on $(X,\gamma)$ define for 
$i\in I$ the $\ZZ$-action  $\tilde e_i^\bullet: \ZZ\times 
X_{\star}(T')\to X_{\star}(T')$ 
by the formula
\begin{equation}
\label{eq:tropicalization}
\tilde e_i^\bullet =Trop(e_{i;\theta}) \ ,
\end{equation}
where $e_{i;\theta}:\GG_m\times T'\to T'$ 
is the positive morphism defined by 
(\ref{eq:positive ei}). 

Also denote $\tilde \gamma_{\theta}:
=Trop(\gamma\circ \theta):X_{\star}(T')\to \Lambda^\vee$. 

Obviously, the map $\tilde e_i^{\,1}:X_{\star}(T')\to X_{\star}(T')$ is a 
bijection for $i\in I$. The bijections 
$\tilde e_i:=\tilde e_i^1$, $i\in I$,
define a free combinatorial pre-crystal $Trop_\theta ({\mathcal X})$ 
on $(X_{\star}(T'),\tilde \gamma_{\theta})$.

\begin{theorem} 
\label{th:positive structure}
Let ${\mathcal X}$ be a geometric crystal, and let 
$\theta$ be a positive structure on $X$. Then 
$Trop_\theta ({\mathcal X})$ is a free combinatorial $W$-crystal. 
\end{theorem}
 
We call this free combinatorial pre-crystal 
$Trop_\theta ({\mathcal X})$ the
{\it tropicalization of ${\mathcal X}$ with respect to the positive
structure $\theta$}.

For  a reduced sequence $\ii=(i_1,\ldots,i_l)$ 
define a morphism $\theta_{\ii}:(\GG_m)^l\to G/B$
by 
\begin{equation}
\label{eq:Kashiwara theta}
\theta_\ii(c_1,\ldots,c_l):=
x_{i_1}(c_1)\cdots x_{i_l}(c_l)\cdot \overline {s_{i_l}}\cdot \overline {s_{i_{l-1}}}\cdots \overline {s_{i_1}}\cdot B \ .
\end{equation}

The following theorem was proved in \cite{bz3} in a 
slightly different form. 

\begin{theorem} 
\label{th:positive structure on X(infty)} For each $\ii\in R(w_0)$ the 
morphism 
$\theta_\ii$ is a positive structure on 
the geometric crystal ${\mathcal X}_{w_0}$, and the 
tropicalization of ${\mathcal X}_{w_0}$ with respect to $\theta_\ii$ 
is equal to the free combinatorial $W$-crystal ${\mathcal B}_{\ii}$. 
All these positive structures
$\theta_\ii$ are equivalent to each other.

\end{theorem} 

\noindent {\bf Remark}.  In \cite{lu} the Lusztig introduced morphisms 
$\theta^\ii: (\GG_m)^{l(w_0)}\to G/B$:
\begin{equation}
\label{eq:Lusztig theta}
\theta^\ii(c_1,\ldots,c_l)
=y_{i_1}(c_1)\cdots  y_{i_l}(c_l)\cdot B \ ,
\end{equation}
and proved 
that these morphisms are related to each other in the same way as 
the corresponding parametrizations of  the canonical basis for 
$U_q({}^L\gg)$. We discuss similar morphisms $\pi^\ii$ in section 
\ref{subsect:positive unipotent structures}. 

It was shown in \cite{bz3} that the morphisms 
$\theta^\ii:(\GG_m)^{l(w_0)}\to G/B$, 
$\ii\in R(w_0)$ are also positive structures on 
${\mathcal X}_{w_0}$, and, moreover, these positive structures are
equivalent to  $\theta_\ii,\ii\in R(w_0)$.

\subsection {Trivialization of geometric crystals} 

\label{subsect:trivialization}

Let  ${\mathcal X}$ be a geometric $G$-crystal built on $(X,\gamma)$. 
Without loss of generality we may assume that $\gamma$ is a regular 
surjective morphism $X\to T$.  
Denote by $X_0=\gamma^{-1}(e)$ the fiber over the unit  
$e\in T$. 

By definition of the action of $W$ on $X$ (see 
Proposition \ref{pr:W-action}), 
$w(x_0)=x_0$  for $x_0\in X_0$, $w\in W$. 

\noindent {\bf Definition}.
A {\it trivialization} of a geometric crystal  ${\mathcal X}$ 
is a $W$-invariant rational projection
$\tau:X\to X_0$.

It is easy to see that the following formula defines a trivialization 
for any geometric $SL_2$-crystal ${\mathcal X}$:
$$\tau(x)=e_i^{\frac{1}{\omega_1(\gamma(x))}}(x)$$
where $I=\{1\}$ and $\omega_1$ is the only fundamental weight of $T$.

In the case when $G=SL_3$  we can construct two different 
trivializations $\tau,\tau':X\to X_0$
for any geometric $SL_3$-crystal ${\mathcal X}$. 
These trivializations are given by 
$$\tau(x)=e_1^{
\frac{1}{\omega_2(\gamma(x))}}
e_2^{\frac{1}{\omega_2(\gamma(x))}}
e_1^{\frac{\omega_2(\gamma(x))}{\omega_1(\gamma(x))}}(x) \ ,$$
$$\tau'(x)= e_2^{\frac{1}{\omega_1(\gamma(x))}}
e_1^{\frac{1}{\omega_1(\gamma(x))}}
e_2^{
\frac{\omega_1(\gamma(x))}{\omega_2(\gamma(x))}}(x) \ ,$$ 
where $\omega_1,\omega_2\in \Lambda$ are the fundamental weights. 
 
{\bf Warning}. The morphism $X\to X_0$: 
$x\mapsto e_1^{\frac{1}{\omega_1(\gamma(x))}}
e_2^{\frac{1}{\omega_2(\gamma(x))}}(x)$ 
is not $W$-invariant. 

\bigskip

Actually for any $r>1$ and any  geometric $SL_{r+1}$-crystal 
${\mathcal X}$ we can construct two trivializations 
$\tau,\tau':X\to X_0$ of ${\mathcal X}$ in the following way. 
Let $\omega_1,\ldots,\omega_r\in \Lambda$ be the fundamental weights 
ordered in the standard way and $\omega_i(x)=\omega_i(\gamma(x))$ 
for $i=1,\ldots,r$. 

For every geometric $SL_{r+1}$-crystal ${\mathcal X}$
define a morphism $\tau:X\to X_0$ by the formula 
$$\tau(x)=\left(e_1^{\frac{1}{\omega_r(x)}}\cdots 
e_r^{\frac{1}{\omega_r(x)}}\right)
\left(e_1^{\frac{\omega_r(x)}{\omega_{r-1}(x)}}\cdots 
e_{r-1}^{\frac{\omega_r(x)}{\omega_{r-1}(x)}}\right) 
\cdots \left(e_1^{\frac{\omega_3(x)}{\omega_{2}(x)}}
e_2^{\frac{\omega_3(x)}{\omega_{2}(x)}}\right)
\left(e_1^{\frac{\omega_2(x)}{\omega_1(x)}}\right)(x) \ .$$

\begin{theorem} 
\label{th:trivialization tau}
This morphism $\tau:X\to X_0$ 
is a trivialization of ${\mathcal X}$. 
\end{theorem}

The formula for the second trivialization $\tau':X\to X_0$ is 
obtained from $\tau$ by applying the automorphism of 
the Dynkin diagram  which exchanges $i$ with $r+1-i$.

\smallskip

We do not know whether trivializations of geometric $G$-crystals 
exist for other reductive groups. We expect that for the 
{\it unipotent crystals}  one  can find a trivialization.

\section{Unipotent crystals}

\label{sect:unipotent crystals}

\subsection{Definition of unipotent crystals and their product}

\label{subsect:definition unipotent}
As we have seen, geometric crystals are geometric 
analogues of free combinatorial $W$-crystals.
Our next task is to introduce a  geometric analogue of 
crystal bases. 

\noindent {\bf Definition}. A $U$-{\it variety} ${\bf X}$  is a pair 
$(X,\alpha)$, where $X\in Ob({\mathcal V})$ and 
$\alpha:U\times X\to X$ is a rational unital $U$-action on $X$ such that 
the restriction of $\alpha$ to each $U_i\times X$, $i\in I$, is a rational action 
$U_i\times X\to X$. 

For $U$-varieties ${\bf X},{\bf Y}$ we say that a rational morphism 
${\bf f}:X\to Y$ is  $U$-{\it morphism} if it commutes with the $U$-actions.

It is well-known that the multiplication in $G$ induces a 
birational isomorphism $B^-\times U\widetilde \to G$. 
Denote by ${\bf g}$ the inverse birational isomorphism:  
\begin{equation}
\label{eq:gauss}
{\bf g}:G\widetilde \to  B^-\times U \ .
\end{equation}

Let $\pi^-:G\to B^-$ and $\pi:G\to U$ be the rational morphisms defined by 
 $$\pi^-=pr_1\circ {\bf g},\pi=pr_2\circ {\bf g} \ .$$
By definition, $dom(\pi^-)=dom(\pi)=B^-\cdot U$.

Passing to the quotient, we obtain a birational isomorphism $G/U\widetilde \to B^-$. 
Therefore, the natural left action of $U$ on
$G/U$ defines a left $U$-action $\alpha_{B^-}:U\times B^-\to B^-$. This action satisfies
for $u\in U, b\in B^-$:
\begin{equation}
\label{eq:the action on B-}
\alpha_{B^-}(u,b)=\pi^-(u\cdot b)=u\cdot b \cdot (\pi(u\cdot b))^{-1} \ .
\end{equation}

In particular, the pair ${\bf B}^-:=(B^-,\alpha_{B^-})$ is a $U$-variety. 

\begin{lemma} 
\label{le:elementary pi}
For $u\in U^-$, $t\in T$ one has:
\begin{equation}
\label{eq:elementary pi}
\pi(x_i(a)\cdot  u\cdot t)
=x_i\left((a^{-1}+\chi^-_i(u))^{-1}\cdot \alpha_i(t^{-1})\right)
\end{equation}
\end{lemma}

\begin{lemma} 
\label{le:action on B-} 
$~~~$

\noindent (a) For $b=u\cdot t$, $u\in U^-, t\in T$,  
$a\in \GG_a$ and $i\in I$ we have:
\begin{equation}
\label{eq:simple generator alpha}
\alpha_{B^-}(x_i(a),b)
=x_i(a)\cdot b\cdot x_i\left(-\frac{a}{(1+a\chi^-_i(u))\alpha_i(t)}
\right)\  .
\end{equation}

\noindent (b) Every $U$-orbit in ${\bf B}^-$ is the 
intersection of $B^-$ with a 
$U\times U$-orbit in $G$.  

\end{lemma}

\noindent {\bf Definition}.  A {\it unipotent crystal} is a 
pair $({\bf X},{\bf f})$,  
where ${\bf X}$ is a $U$-variety and ${\bf f}:{\bf X}\to {\bf B^-}$ 
is a $U$-morphism.

We denote by  $U-{\mathcal Cryst}$ the category whose 
objects are unipotent $G$-crystals and arrows are dominant 
rational morphisms. 

For any $({\bf X},{\bf f}_X),({\bf Y},{\bf f}_Y)\in Ob(U-{\mathcal Cryst})$ 
define a rational morphism $\alpha:U\times X\times Y\to X\times Y$ 
by the formula:
\begin{equation}
\label{eq:unipotent group action on the product}
\displaystyle{\alpha(u,(x,y))
:=\left(u(x),\left(\pi(u\cdot{\bf f}_X(x))\right)(y)\right)} \ ,
\end{equation}

where $u(x)=\alpha(u,x)$. We will often write $u(x,y)$ instead of $\alpha(u,(x,y))$

\begin{theorem}
\label{th:product}

$~~~$

\noindent (a) The morphism 
$\alpha:U\times X\times Y\to X\times Y$ defined above is a rational 
$U$-action on $X\times Y$. 

\noindent (b)  Let ${\bf m}:B^-\times B^-\to B^-$ be the 
multiplication morphism. Let ${\bf f}={\bf f}_{X\times Y}:X\times Y\to B^-$ 
be the rational morphism defined by 
${\bf f}={\bf m}\circ ({\bf f}_X \times {\bf f}_Y)$.
Then  ${\bf f}_{X\times Y}$ is a $U$-morphism.

\end{theorem}

We denote the $U$-variety  $(X\times Y,\alpha_{X\times Y})$ by 
${\bf X}\times_{\bf f} {\bf Y}$. According to the Theorem 
the pair $({\bf X}\times_{\bf f} {\bf Y},{\bf f}_{X\times Y})$  
is a unipotent crystal. We call it the {\it product} of 
$({\bf X},{\bf f}_X)$ and 
$({\bf Y},{\bf f}_Y)$ and  denote it by  
$({\bf X},{\bf f}_X)\times ({\bf Y},{\bf f}_Y)$.

\smallskip

\noindent {\bf Remark}. Product of unipotent crystals is 
analogous to the tensor product of Kashiwara's crystals defined in 
\cite {k93}. This analogy is made precise in sections 
\ref{subsect:unipotent to geometric} and \ref{subsect:positive unipotent}
below.
 
\begin{proposition} 
\label{pr:monoidal category}
Product of unipotent crystals is associative. 
\end{proposition}

\noindent {\bf Remark}. The above results define a strict monoidal 
structure on the category $U-{\mathcal Cryst}$. 

\smallskip

The product of unipotent crystals is not commutative in general. 
For any family $({\bf X}_k,{\bf f}_k)\in Ob(U-{\mathcal Cryst})$, 
$k=1,\ldots,l$ we denote by $\prod\limits_{k=1}^l ({\bf X}_k,{\bf f}_k)$ 
the product 
$({\bf X}_1,{\bf f}_1)\times \cdots \times ({\bf X}_l,{\bf f}_l)$. 

\medskip

\noindent {\bf Example}. Denote by ${\bf U}$ 
the pair $(U,\alpha_U)$, where $\alpha_U:U\times U\to G$ 
is the left action.
Let $pr_e^U:U\to B^-$ be the projection  on  the unit  $e$. 
Clearly, $({\bf U},pr_e^{U})\in Ob(U-{\mathcal Cryst})$.  
Similarly, denote by ${\bf G}$ the pair $(G,\alpha_G)$, 
where $\alpha_G:U\times G\to U$ is the left action. 
Clearly, $({\bf G},\pi^-)\in Ob(U-{\mathcal Cryst})$.

\begin{proposition} 
\label{pr:G is crystal}
The birational isomorphism ${\bf g}$ induces the 
isomorphism in $U-{\mathcal Cryst}$:
$$({\bf G},\pi^-)\widetilde \to ({\bf B}^-,{\rm id}_{B^-})\times ({\bf U},pr_e^{U}) \ .$$
More generally, taking the right quotient by any subgroup $U'\subset U$, 
we obtain the birational isomorphism 
${\bf g}_{U'}:G/U'\widetilde \to B^-\times U/U'$ which induces the isomorphism 
in $U-{\mathcal Cryst}$: 
$$({\bf G/U'},\pi^-_{U'})\widetilde \to 
({\bf B}^-,{\rm id}_{B^-})\times ({\bf U/U'},pr_e^{U/U'})$$
where $\pi^-_{U'}=pr_1\circ {\bf g}_{U'}$. 
\end{proposition}

The following result shows that the unipotent crystal $({\bf G},\pi^-)$ is universal.

\begin{lemma} 
\label{le:universal G}
For every $U$-equivariant rational morphism 
${\bf f}_U:U\to B^-$ there is an element $\tilde w\in Norm_G(T)$ such that 
$${\bf f}_U(u)=\pi^-(u\cdot \tilde w)$$
for every $u\in U$. 
\end{lemma}

\subsection{From unipotent $G$-crystals to unipotent $L$-crystals}

\label{subsect:L-crystals}

Throughout this section we fix a subset $J$ of $I$.  
Let $L=L_J$ be the  Levi subgroup of $G$ 
generated by $T$ and by $U_j,U^-_j$, $j\in J$. 

Let $P:=L\cdot U$ and $P^-:=U^-\cdot L$.  By definition, 
$P$ and $P^-$ are 
parabolic subgroups of $G$ such that  
$P\supset B$, $P^-\supset B^-$, and 
$P\cap P^-=L$.

Let $U_P$ and $U^-_P$ be the respectively the unipotent radicals of 
$P$ and $P^-$. 
Let $U_L=L\cap U$, $U_L^-=L\cap U^-$. Then $U_L$ and $U^-_L$ are 
the opposite unipotent radicals of $L$. 

Denote $B^-_L:=B^-\cap L$, and $B_L=B\cap L$. The open inclusion 
$B^-_L\hookrightarrow L/U_L$ induces a rational action of $U_L$ on 
$B^-_L$ which we denote by $\alpha_L:U_L\times B^-_L\to B^-_L$. 
Let ${\bf p}^-={\bf p}^-_L:B^-\to B^-_L$ be the canonical projection.  By definition, 
${\bf p}_L^-$ commutes with the rational action of $U_L$.  
In particular, $({\bf B}^-,{\bf p}_L^-)\in Ob(U_L-{\mathcal Cryst})$. 

For any $U$-variety ${\bf X}$ we denote by ${\bf X}|_L$ the $U_L$-variety 
obtained by the restriction of the $U$-action 
$\alpha:U\times X\to X$ to $U_L$. 

\begin{lemma} The mapping 
$({\bf X},{\bf f}_X)\mapsto ({\bf X},{\bf f}_X)|_L:=
({\bf X}|_L,{\bf p}_L^-\circ {\bf f}_X)$ defines a functor of monoidal 
categories $|_L:U-{\mathcal Cryst}\longrightarrow U_L-{\mathcal Cryst}$. 

\end{lemma}  

We call the unipotent $L$-crystal 
$({\bf X},{\bf f}_X)|_L$ the {\it restriction} of $({\bf X},{\bf f}_X)$ to $L$.

\subsection{From unipotent $G$-crystals to  geometric $L$-crystals}

\label{subsect:unipotent to geometric}

For each  unipotent $G$-crystal $({\bf X},{\bf f})$ define a morphism 
$\gamma=\gamma_X:X\to T$ by $\gamma=pr_T\circ {\bf f}$.  
For each $i\in I$ define the function  
$\varphi_i=\varphi_i^X$ by 
\begin{equation}
\label{eq:phi on X}
\varphi_i:=\overline \chi_i^{\,-}\circ {\bf f}_X \ . 
\end{equation}
 
Let ${\rm supp}({\bf X},{\bf f})$ be set of all 
those $i\in I$ for which $\varphi_i^X\not \equiv 0$. 
We call this set the {\it support} of the unipotent crystal 
$({\bf X},{\bf f}_X)$.
For $i\in {\rm supp}({\bf X},{\bf f})$ 
define the morphism 
$e_i:\GG_m\times X\to X$  by
\begin{equation}
\label{eq:action of ei2}
e_i^c(x)= x_i\left(\frac{c-1}{\varphi_i(x)}\right)(x) \ .
\end{equation}

It is easy to see that each $e_i$ is a rational action of $\GG_m$ on $X$.

\begin{theorem} 
\label{th:the first functor}
For any $({\bf X},{\bf f})\in Ob(U-{\mathcal Cryst})$ 
the actions
$e_i:\GG_m\times X\to X$, $i\in {\rm supp}({\bf X},{\bf f})$ define a 
geometric $L_J$-crystal on
$(X,\gamma_X)$, where $J={\rm supp}({\bf X},{\bf f})$.   
\end{theorem}  

We denote this geometric $L_J$-crystal by ${\mathcal X}_{\rm ind}$ 
and call it the {\it geometric crystal induced by}  
$({\bf X},{\bf f})$.

\smallskip

\noindent {\bf Remark}. For the unipotent crystal 
$(G/B,[{\bf j}_0]^{-1})$ (defined in section \ref{subsect:Xinfty}) 
Theorem \ref{th:the first functor}, specializes to Theorem
\ref{th:X(infinity)}.  

\smallskip

\noindent {\bf Examples}.

\noindent 1.  For the unipotent $G$-crystal $({\bf B}^-,{\rm id}_{B^-})$ 
the actions  $e_i:\GG_m \times B^-\to B^-$, $i\in I$, are given by  
\begin{equation}
\label{eq:simple multiplicative generator ei}
e_i^c(b)=x_i\left(\frac{c-1}{\varphi_i(b)}\right)\cdot 
b\cdot x_i\left(\frac{c^{-1}-1}
{\varphi_i(b)\alpha_i(\gamma(b))}\right) \ .
\end{equation}

\noindent 2. For the unipotent crystal $({\bf G},\pi^-)$ 
the actions  $e_i:\GG_m \times G\to G$, $i\in I$, are given by  
\begin{equation}
\label{eq:simple multiplicative generator ei on G}
e_i^c(g)=x_i\left(\frac{c-1}{\overline \chi_i(g)}\right)\cdot g \ .
\end{equation}

\medskip

\begin{lemma} 
\label{le:geometric crystal on the product}
For  $({\bf X},{\bf f}_X),({\bf Y},{\bf f}_Y)\in Ob(U-{\mathcal Cryst})$  
put $({\bf Z},{\bf f}_Z):=({\bf X},{\bf f}_X)\times ({\bf Y},{\bf f}_Y)$, 
where $Z=X\times Y$, and let 
${\mathcal Z}_{ind}$ be the geometric crystal 
on $(Z,\gamma_Z)$ induced by $({\bf Z},{\bf f}_Z)$. We have:

\noindent (a) $\gamma_Z={\bf m}\circ (\gamma_X\circ \gamma_Y)$, where 
${\bf m}:T\times T\to T$ is the multiplication morphism. 

\noindent (b) For each $i\in I$, $(x,y)\in Z$: 
\begin{equation}
\label{eq:phi of the product}
\varphi_i^Z(x,y)=\varphi_i^X(x)
+\frac{\varphi_i^Y(y)}{\alpha_i(\gamma_X(x))}
\end{equation}
which implies that ${\rm supp}({\bf Z},{\bf f}_Z)
={\rm supp}({\bf X},{\bf f}_X)\cup {\rm supp}({\bf Y},{\bf f}_Y)$.

\noindent (c) For any $i\in {\rm supp}({\bf Z},{\bf f}_Z)$ the action 
$e_i:\GG_m\times Z\to Z$ is given by the formula:
$e_i^c(x,y)=(e_i^{c_1}(x),e_i^{c_2}(y))$, 
where
\begin{equation}
\label{eq:tensor product of crystals}
c_1=\frac{c\varphi_i(x)\alpha_i(\gamma(x))+\varphi_i(y)}
{\varphi_i(x)\alpha_i(\gamma(x))+\varphi_i(y)},
~~c_2=\frac{\varphi_i(x)\alpha_i(\gamma(x))+\varphi_i(y)}
{\varphi_i(x)\alpha_i(\gamma(x))+c^{-1}\varphi_i(y)} \ .
\end{equation}
\end{lemma}

\noindent {\bf Remark}.  The formula 
(\ref{eq:tensor product of crystals}) for the action of 
$e_i$ on $X\times Y$ 
is analogous to the formula in \cite{k93}, section 1.3,  
for the tensor product of Kashiwara's crystals.

\subsection{Positive unipotent crystals and duality}

\label{subsect:positive unipotent}

Let $({\bf X},{\bf f})\in Ob(U-{\mathcal Cryst})$, and let $T'$ be an 
algebraic torus of the same dimension as $X$. A birational isomorphism 
$\theta:T'\widetilde \to X$ is called 
a {\it positive structure} on $({\bf X},{\bf f})$ if the following 
two conditions are satisfied. 

\noindent (1) The isomorphism $\theta$  is a positive structure on  the 
induced geometric $L_J$-crystal ${\mathcal X}_{\rm ind}$, where 
$J={\rm supp}({\bf X},{\bf f})$. 

\noindent (2) For any $i\in J$ the function 
$\varphi^X_i\circ \theta$ on $T'$ is 
positive. 

\begin{theorem} 
\label{th:positive unipotent}
Let $({\bf X},{\bf f}_X), ({\bf Y},{\bf f}_Y)\in 
Ob(U-{\mathcal Cryst})$  and $\theta_X:T'\to X,~\theta_Y:T''\to Y$ 
be respectively the positive structures. 
Then the birational isomorphism $\theta_{X\times Y}:=
\theta_X\times \theta_Y$ 
is a positive structure on the product 
$({\bf X},{\bf f}_X)\times ({\bf Y},{\bf f}_Y)$. 

\end{theorem}

For a  geometric pre-crystal ${\mathcal X}$ on $(X,\gamma)$ denote by
$\gamma^* :X\to T$ the morphism  $\gamma^*(x)=(\gamma(x))^{-1}$ and 
consider the {\it dual} geometric pre-crystal ${\mathcal X}^*$ on
$(X,\gamma^*)$  by defining 
$\displaystyle{(e_i^*)^c(x)=e_i^{c^{-1}}(x)}$. 

\medskip

Given a geometric pre-crystal ${\mathcal X}$ on $(X,\gamma)$, 
for any morphism $\theta:T'\to X$ define a morphism
$\theta^*:T'\to X$ as the composition of $\theta$ with the inverse 
${}^{-1}:T'\to T'$. 

The following fact is obvious. 

\begin{lemma} For any geometric crystal ${\mathcal X}$ 
the dual geometric pre-crystal ${\mathcal X}^*$ is a geometric crystal.
For any positive structure $\theta$ on ${\mathcal X}$ the 
morphism $\theta^*$ is a positive structure on ${\mathcal X}^*$. 
\end{lemma}

\noindent {\bf Remark}. Duality in  geometric crystals is an  
analogue of duality in Kashiwara's crystals (see \cite{k93}  
and the Appendix below). 
More precisely, the tropicalization of ${\mathcal X}^*$ with respect to 
the positive  structure $\theta^*$ is the free combinatorial $W$-crystal 
dual to the tropicalization ${\mathcal X}$ with respect to $\theta$. 

\smallskip

Let $({\bf X},{\bf f})$ be a unipotent crystal. 
Define a morphism ${\bf f}^*:X\to B$ by ${\bf f}^*(x)=({\bf f}(x))^{-1}$ 
and a rational morphism
$\alpha^*:U\times X\to X$ by  
$$\alpha^*(u,x):=\pi\left(u\cdot{\bf f}^*(x))\right)(x) \ .$$

\begin{proposition}  
\label{pr:dual unipotent crystal}
The morphism $\alpha^*$ is a rational action of $U$ on $X$ and 
${\bf f}^*$ is a $U$-morphism with respect to this action. 
\end{proposition}

Denote the pair $(X,\alpha^*)$ by ${\bf X}^*$. 
Then the pair $({\bf X}^*,{\bf f}^*)$ is a  unipotent crystal. 
We call it the {\it dual unipotent crystal}  
of $({\bf X},{\bf f})$ and denote it by $({\bf X},{\bf f})^*$.

\noindent {\bf Remark}. It follows from the definition that the mapping  
$({\bf X},{\bf f})\mapsto ({\bf X},{\bf f})^*$ is an involutive functor 
${}^*:U-{\mathcal Cryst}\longrightarrow U-{\mathcal Cryst}$.  

\smallskip

\begin{theorem} 
\label{th:duality of product}
For any unipotent crystals  $({\bf X},{\bf f}_X), 
({\bf Y},{\bf f}_Y)$ 
the permutation of the factors 
$(12):X\times Y\to Y\times X$
induces an isomorphism of unipotent crystals:
\begin{equation}
\left(({\bf X},{\bf f}_X)\times({\bf Y},{\bf f}_Y)\right)^*\widetilde \to 
({\bf Y},{\bf f}_Y)^*\times ({\bf X},{\bf f}_X)^*  \ .
\end{equation}
\end{theorem}

\noindent {\bf Remark}. This theorem implies that the functor 
${}^*:U-{\mathcal Cryst}\longrightarrow U-{\mathcal Cryst}$ reverses the 
monoidal structure on $U-{\mathcal Cryst}$. 

\smallskip

\begin{theorem} 
\label{th:dual unipotent to geometric} 
Let $({\bf X},{\bf f})$ be a unipotent crystal and 
${\mathcal X}$ the 
induced geometric crystal. Then the geometric crystal induced by
$({\bf X},{\bf f})^*$ is equal to ${\mathcal X}^*$. 
\end{theorem}

\subsection {Diagonalization of products of unipotent crystals}

\label{subsect:diagonalization}

Let ${\bf v}$  be the regular morphism $Bw_0B\to U$ defined by
${\bf v}(ut\overline{w_0}u')=uu'$
for any $u,u'\in U$, $t\in T$. 

By definition, ${\bf v}$ is two-sided $U$-equivariant.

\noindent {\bf Remark}. For $G=SL_2$ the map ${\bf v}:G\to U$ is given by 
$${\bf v}\begin{pmatrix} a & b \\ c & d \end{pmatrix}=
\begin{pmatrix} 1 & \frac{a+d}{c} \\ 0 & 1 \end{pmatrix} $$

\smallskip

We call a unipotent crystal $({\bf X},{\bf f}_X)$ {\it non-degenerate} 
if $ran({\bf f}_X)$ intersects $Bw_0B$ non-trivially. 

For any non-degenerate unipotent crystal 
$({\bf X},{\bf f}_X)\in Ob(U-{\mathcal Cryst})$ and 
any $({\bf Y}, {\bf f}_Y)\in Ob(U-{\mathcal Cryst})$ 
define the rational morphism
$F=F_{X\times Y}:X\times Y\to X\times Y$ by 
$F(x,y)=(x,v_x(y))$, 
where $v_x:={\bf v}({\bf f}_X(x))$. 
Clearly, $F$ is a birational isomorphism  
$X\times Y\widetilde \to X\times Y$, and its 
inverse is given by $F^{-1}(x,y)=
(x,(v_x)^{-1}(y))$.

Denote by $\delta=\delta_{X\times Y}:
U\times X\times Y\to X\times Y$ the diagonal action of $U$:
$$(u,(x,y))\mapsto (u(x),u(y)) \ .$$

\begin{proposition} 
\label{pr:diagonalization} Let  
$({\bf X},{\bf f})$ be a non-degenerate unipotent crystal, 
and $({\bf Y}, {\bf f}_Y)$ be any unipotent crystal. Then
\begin{equation}
\label{eq:diagonalization}
F\circ \alpha=\delta\circ ({\rm id}_U\times F) \ . 
\end{equation} 
In other words  the action 
$\alpha=\alpha_{X\times Y}:U\times X\times Y\to X\times Y$  
is diagonalized by $F$.
\end{proposition}

\subsection {Unipotent action of the Weyl group}

\label{subsect:adjoint model}

For each $J\subset I$ let $W_J$ be the subgroup of $W$ generated by 
$s_j,j\in J$. By definition, $W_J$ is the Weyl group of the 
standard Levi subgroup  $L=L_J$. So we sometimes denote $W_J$ by $W_L$.

\begin{lemma} 
\label{pr:induce unipotent action}
Let $({\bf X},{\bf f})$ be a unipotent $G$-crystal, 
$J=supp ({\bf X},{\bf f})$ and ${\mathcal X}_{\rm ind}$ be the 
induced geometric $L_J$-crystal. 
Then the formula (\ref{eq:W-action}) defines rational action of $W_J$ 
on $X$. 
\end{lemma}

We call this action of $W_J$ on $X$ the {\it unipotent action} of $W_J$ 
on $({\bf X},{\bf f}_X)$ (or, simply, {\it unipotent action} of $W_J$ on $X$).

Recall from section \ref{subsect:definition unipotent}
that $({\bf B}^-,{\rm id}_{B^-})$ is a unipotent crystal.

\begin{lemma} 
The unipotent action of $W$ on 
$B^-$ satisfies for $i\in I$: 
\begin{equation}
\label{eq:simple generator si}
s_i(b)=x_i(a)\cdot b\cdot (x_i(a))^{-1} \ ,
\end{equation}
where $a=\frac{1-\alpha_i(\gamma(b))}{\varphi_i(b)\alpha_i(\gamma(b))}$. 
\end{lemma}

\noindent {\bf Remark}. This lemma implies that $dom(s_i)=
\{b\in B:\overline \chi_i^{\,-}(b)\ne 0\}$. 

\smallskip

For each $w\in W$ let ${\rm supp}(w)$ be the minimal subset $J\subset I$ such 
that $w\in W_J$. For example, ${\rm supp}(s_i)=\{i\}$.

\begin{proposition} 
\label{pr:uw}
For each $w\in W$ there is a rational morphism  
${\bf u}_w: B^-\to U$  such that: 

\noindent (a)  
For any $w'\in W$ such that ${\rm supp}(w')\supset {\rm supp}(w)$ each rational 
$U$-orbit of the form $t\cdot U\overline {w'}U\cap B^-$, $t\in T$ 
intersects  $dom({\bf u}_w)$ non-trivially. 

\noindent (b) The unipotent action of $W$ on $B^-$ is given by 
\begin{equation}
\label{eq:unipotent w}
(w,b)\mapsto w(b)={\bf u}_w(b)\cdot b\cdot ({\bf u}_w(b))^{-1} \ .
\end{equation}

\end{proposition}

For a any unipotent crystal $({\bf X},{\bf f}_X)$ and any 
$w\in W$ with ${\rm supp}(w)\subset {\rm supp}({\bf X},{\bf f}_X)$  
we define a rational morphism ${\bf u}_w^X:X\to U$ by 
${\bf u}_w^X={\bf u}_w\circ {\bf f}_X$.

\begin{theorem} Let  $({\bf X},{\bf f}_X)\in Ob(U-{\mathcal Cryst})$ such that 
${\rm supp}({\bf X},{\bf f}_X)=J$. Then the crystal action of 
$W_J$ on $X$ is given by 
\begin{equation}
\label{eq:the w}
w(x)=\left({\bf u}_w^ X(x)\right)(x)
\end{equation} 
for $w\in W_J$. 
\end{theorem}

\begin{proposition} 
\label{co:w n the product} Let $({\bf X}_k,{\bf f}_{X_k})\in 
Ob(U-{\mathcal Cryst})$, $k=1,\ldots,l$, and 
let $J=\cup_k {\rm supp} ({\bf X}_k,{\bf f}_{X_k})$. 
Then the unipotent action of $W_J$ on $X=X_1\times \cdots  \times X_l$ is given by 
$$w(x_1,\ldots,x_l)= \left(u(x_1),\ldots,u(x_l)\right)$$
for $w\in W$, where 
$u={\bf u}_w^X(x_1 \ldots x_l)$ is as in (\ref{eq:the w}). 

\end{proposition}

\section{Examples of unipotent crystals}

\label{sect:examples}
\subsection{Standard unipotent crystals}

\label{subsect:standard crystals}
For $w\in W$ let ${\mathcal O}(w):=U\overline wU/U$. Clearly, 
${\mathcal O}(w)$ is a left $U$-orbit in $G/U$. Define  
$U(w):=U\cap \overline wU\overline w^{\,-1}$.

\begin{lemma} The mapping $U\to {\mathcal O}(w)$ defined by 
$u\mapsto u\overline wU$ for $u\in U$ is left $U$-equivariant 
and surjective. It induces  the isomorphism of homogeneous $U$-varieties:
$$\tilde \eta^w:U/U(w)\widetilde \to {\mathcal O}(w) \ .$$
\end{lemma}

\noindent {\bf Remark}. Under the natural map $G\to G/B$  each orbit 
${\mathcal O}(w)$ is identified with the corresponding Schubert cell. 
In particular, $\dim~{\mathcal O}(w)=\dim~U/U(w)=l(w)$. 

\smallskip

For $w\in W$ define $U^w:=U\cap B^- w^{\,-1} B^-$
and $B^-_w:=U\overline wU\cap B^-$. 

\noindent {\bf Example}. For $G=SL_2$, $w=w_0$ the sets
$U^w$ and $B^-_w$  consist respectively of the matrices of the form:
$$\begin{pmatrix} 
1 & c \\
0 & 1
\end{pmatrix}, ~~\begin{pmatrix} 
c & 0 \\
1 & c^{-1}
\end{pmatrix} \ ,$$
$c\in \GG_m$.  

\medskip

Denote by ${\bf j}^w:U^w\to U/U(w)$ the morphism induced by the natural 
inclusion $U^w\hookrightarrow U$, and  by ${\bf j}_w:B^-_w\to 
{\mathcal O}(w)$ the 
morphism induced by the natural inclusion $B^-_w\hookrightarrow U\overline w U$.

By definition, the restrictions of $\pi$ and 
$\pi^-$ to $B^-\cdot U$ are the regular projections 
$B^-\cdot U\to U$ and $B^-\cdot U\to
B^-$ respectively. 
Since $B^-\overline w^{\,-1}B^-\cdot 
\overline w\subset B^-\cdot U$ we
can define a regular morphism  
$\eta^w:U^w\to B^-$ by
$$\eta^w(u)=\pi^-(u\overline w)$$
for $u\in U^w$.

\begin{proposition} 
\label{pr:commutative diagram bz3}                             
$~~~$                                   

\noindent (a) The morphisms ${\bf j}_w$ and ${\bf j}^w$ 
are the open inclusions 
$B^-_w\hookrightarrow {\mathcal O}(w)$ and 
$U^w\hookrightarrow U/U(w)$ respectively. 

\noindent (b) The morphism 
$\eta^w$ is a biregular isomorphism $U^w\widetilde \to B^-_w$. 
The inverse isomorphism 
$\eta_w=(\eta^w)^{-1}:B^-_w\widetilde \to U^w$ is given by 
\begin{equation}
\label{eq:inverse eta}
\eta_w(b)=\pi(\overline w\cdot b^{-1})^{-1} \ .
\end{equation}

\noindent (c) The following diagram is commutative:
\begin{equation}
       \begin{CD} 
       U/U(w)   @>\tilde \eta^w>>  {\mathcal O}(w)    \\ 
            @AA {\bf j}^w A   @ AA {\bf j}_w A  \\
        U^w          @>\eta^w>>         B^-_w      
        \end{CD} 
       \label{eq:commutative diagram} ~~~\ .
   \end{equation}        

\end{proposition}

The birational isomorphisms $[{\bf j}^w]$ and 
$[{\bf j}_w]$ define for any $w\in W$ rational unital $U$-actions
$\alpha^w:U\times U^w\to U^w$ and 
$\alpha_w:U\times B^-_w\to B^-$ respectively.   
We denote by ${\bf U}^w:=(U^w,\alpha^w)$ and 
${\bf B}^-_w:=(B_w,\alpha_w)$ the corresponding $U$-varieties. 
 
Note that  $({\bf U}^w,\eta^w)$ 
and $({\bf B}^-_w,{\rm id}_w)$ are unipotent $G$-crystals. 
We call these unipotent crystals {\it standard}. 

\begin{lemma} 
\label{le:the twist U to B-}
For each $w\in W$ the isomorphism $\eta^w$ 
induces the isomorphism of unipotent crystals 
$({\bf U}^w,\eta^w)\widetilde \longrightarrow  ({\bf B}^-_w,{\rm id}_w)$.
\end{lemma}

Note that the $U$-action on ${\bf B}^-_w$ is given by 
(\ref{eq:the action on B-}). 
In order to describe the $U$-action on ${\bf U}^w$ we need to introduce more notation. 

By Proposition \ref{pr:commutative diagram bz3}, the multiplication 
morphism $U\times U\to U$ 
induces the birational isomorphism $U^w\times U(w)\widetilde \to U$. 
Denote by $\pi^w:U\to U(w)$ the composition of the inverse 
of this isomorphism with $pr_2$. 

By definition, the action $U\times U^w\to U^w$: 
$(u,x)\mapsto u(x)$ is given by: 
\begin{equation}
\label{eq:action of U on Uw}
u(x)=u\cdot x \cdot (\pi^w(u\cdot x))^{-1} \ .
\end{equation}

\begin{lemma}
For any $w, w'\in W$ such that $l(ww')=l(w)+l(w')$ 
and $x\in U^{w'}, y\in U^w$, $u\in U$ we have:
$$\pi^{ww'}(u\cdot x\cdot y)=\pi^w(\pi^{w'}(u\cdot x)\cdot y)$$
or, equivalently, 
$u(x\cdot y)=u(x)\left(\pi^{w'}(u\cdot x)(y)\right)$.

\end{lemma}

Next, we compute the action (\ref{eq:action of U on Uw}) 
explicitly in terms of the isomorphism $\eta^w$. 

\begin{proposition} 

\label{pr:U-action on Uw}
 For each For each $w\in W$,  $u\in U$, and $x\in U^w$ we have:
\begin{equation}
\label{eq:pi for U-action on Uw}
\pi^w(u\cdot x)=\pi\left(Ad~\overline w(\pi(u\cdot b))\cdot 
\pi^-(\overline w\cdot b^{-1})\right)
\end{equation}
where $b=\eta^w(x)$. 

\end{proposition}

\subsection{Multiplication of standard unipotent crystals}

\label{subsect:factorization of unipotent $L$-crystals}

Define a new associative multiplication ${\star}$ on the set $W$
by 
$$w{\star}w'=ww'$$ 
if $l(ww')=l(w)+l(w)$ for $w, w'\in W$ 
and $s_i{\star}s_i=s_i$ for all $i$. 

Under this new operation $W$ is identified with the standard 
multiplicative monoid  of the degenerate Hecke algebra of $W$. 
We denote this monoid by $(W,{\star})$. 

\begin{lemma} 
\label{le:product unipotent bruhat}
For any  $w,w'\in W$  
the multiplication morphism $G\times G\to G$ 
induces a dominant rational morphism
\begin{equation}
\label{eq:product degenerate Bruhat}
\pi^{w',w}:U^{w'}\times U^{w} \to U^{w{\star}w'}  \ . 
\end{equation}
This is a birational isomorphism if and only if $w{\star}w'=ww'$. 

\end{lemma}

It is well-known (see e.g., \cite{fz}, Theorems 1.2, 1.3) 
that for any  $w,w'\in W$  such that $w{\star}w'=ww'$ 
the multiplication morphism $G\times G\to G$ 
induces the open inclusion
\begin{equation}
\label{eq:product Bruhat}
B^-_w\times B^-_{w'}  \to B^-_{ww'} \ .
\end{equation}

The following result deals with a generalization (\ref{eq:product Bruhat}) 
of to any pair $w,w'$.

\begin{proposition} 
 \label{pr:general product Bruhat} 
For any $w,w'\in W$:

\noindent (a) The intersection of $B^-_w\cdot B^-_{w'}$ with 
$B^-_{w{\star}w'}$ 
is a dense subset of $B^-_{w{\star}w'}$. 

\noindent (b) There exists an algebraic sub-torus 
$\widetilde T_{w,w'}\subset T$ such that the restriction of the 
multiplication morphism $G\times G\to G$ 
to $B^-_w\times B^-_{w'}$ is a dominant rational morphism
\begin{equation}
\label{eq:general product Bruhat}
B^-_w\times B^-_{w'} \to B^-_{w{\star}w'}\cdot \widetilde T_{w,w'} \ .
\end{equation}
This is a birational isomorphism if and only if 
$l(w{\star}w')=l(w)+l(w')-dim~\widetilde T_{w,w'}$.

\noindent (c) $\widetilde T_{w,w'}=\{e\}$ if and only if $w\star w=ww'$. 

\noindent (d) $\widetilde T_{s_i,s_i}=\alpha_i^\vee(\GG_m)$ for each $i\in I$. 

\noindent (e) For any $w,w',w''\in W$ one has:
$\widetilde T_{w,w'w''}\cdot \widetilde T_{w',w''}=
\widetilde T_{w{\star}w',w''}\cdot (w'')^{-1}(\widetilde T_{w,w'})$.

\end{proposition}

\begin{corollary} 
\label{co:bruhat in B-}
For any $w,w'\in W$ the morphism (\ref{eq:general product Bruhat}) 
induces the morphism in $U-{\mathcal Cryst}$: $({\bf B}^-_w,{\rm id}_w)\times 
({\bf B}^-_{w'},{\rm id}_{w'})\to 
\left({\bf B}^-_{w{\star}w'}\cdot \widetilde T_{w,w'},
{\rm id}_{B^-_{w{\star}w'}\cdot \widetilde T_{w,w'}}\right)$.
\end{corollary}

Next, we compute a lower bound for each $\widetilde T_{w,w'}$.

For each $w\in W$ let us consider 
a homomorphism of tori $T\to T$ defined by $t\mapsto w(t)\cdot t^{-1}$. 
Denote by $T_w$ the image of this homomorphism. Clearly, $T_w$ 
is a sub-torus of $T$ such that 
$$X_\star(T_w)=(\Lambda^w)^\perp\bigcap \oplus_i \ZZ\alpha_i^\vee \ ,$$ 
where $\Lambda^w=\{\lambda\in \Lambda:w(\lambda)=\lambda\}$. This implies that 
$X_\star(T_w)$ has a $\ZZ$-basis of certain (not necessarily simple) coroots. 

\begin{proposition} 
\label{pr:torus in product of reduced cells} 
For any $w, w'\in W$ we have
$\widetilde T_{w,w'}\supset T_{(w\star w')^{-1}\cdot w\cdot w'}$.

\end{proposition}

{}From now on we will freely use the notation of section 
\ref{subsect:L-crystals}. For any standard Levi subgroups $L\subset L'\subset G$ define 
the elements $w_{L',L},w_{L,L'}\in W$ by 
\begin{equation}
\label{eq:wll'}
w_{L',L}:=w_0^{L'}w_0^L,  ~w_{L,L'}:=w_0^Lw_0^{L'}
\end{equation}
Clearly, $w_{L',L}\cdot w_{L',L}=e$ and 
$l(w_{L',L})=l(w_{L,L'})=l(w_0^{L'})-l(w_0^L)$. 

Note that $U_{L'}(w_{L,L'})=U_L$ and $U_{L'}(w_{L',L})=Ad~\overline {w_{L',L}}(U_L)$.

It turns out that for $w=w_{L',L}$ 
the formula (\ref{eq:pi for U-action on Uw}) simplifies 
when $u\in U_L$. 

\begin{proposition} 
\label{pr:UL-action} Let $L\subset L'$ be  standard Levi subgroups of $G$. Then 
for any $u\in U_L$, $x\in U^{w_{L',L}}$ we have:
$$\pi^{w_{L',L}}(u\cdot x)=Ad~\overline {w_{L',L}}( 
\pi(u\cdot \eta^{w_{L',L}}(x))) \ . $$
\end{proposition}

For each $w\in W$ let $I(w)$ be the set of all 
$i\in I$ such that  $w(\alpha_i)=\alpha_{i'}$ for some $i'\in I$. 
Note that an element $w'\in W$ belongs to  $W_{I(w)}$ if and only if 
$Ad~\overline w (U^{w'})=U^{ww' w^{-1}}$. 

\begin{lemma} 
\label{le:Ad w and pi w}
Let $w\in W$, $w',w''\in W_{I(w)}$. Then
$$Ad~\overline w\circ \pi^{w'}|_{U^{w''}}=\pi^{ww'w^{-1}}|_{U^{ww''w^{-1}}}\circ 
Ad~\overline w \ .$$
\end{lemma}

\begin{theorem} 
\label{th:factorization U} Let $L\subset L'$ be  standard Levi subgroups of $G$, and
let $w \in W$ be any element such that $\overline {w_{L',L}}^{\,-1}\cdot \overline w$ 
centralizes $U_L$. 
Then for any $w'\in W_{I(w)}$  we have:

\noindent (a) The morphism 
$f_{w,w'}:U^{w_{L',L}}\times U^{w'}\to U$ given by 
$(x,y)\mapsto x\cdot 
Ad~\overline w(y)$
is a dominant rational morphism
$U^{w_{L',L}}\times U^{w'}\to U^{(ww'w ^{-1})\star w_{L',L}}$ \ .
This is a birational isomorphism if and only if 
$(ww'w ^{-1})\star w_{L',L}=(ww'w ^{-1})\cdot w_{L',L}$.

\noindent (b) If  $f_{w,w'}$ is a birational isomorphism then it induces
an isomorphism in the category $U_L-{\mathcal Cryst}$:
$$({\bf U}^{w_{L',L}},\eta^{w_{L',L}})|_L\times 
({\bf U}^{w'},\eta^{w'})|_L\widetilde \to 
\left({\bf U}^{(ww'w^{-1})\star w_{L',L}},\eta^{(w w'w ^{-1})\star w_{L',L}}\right)|_L$$
\end{theorem}

\smallskip

Let $G=GL_{m+n}$, $m>0,n>0$. We use the standard labeling of the 
Dynkin diagram of type $A_{m+n-1}$: 
$I=\{1,2,\ldots,m+n-1\}$. Let $J=I\setminus \{m\}$. Clearly, 
$L_J=L_{m,n}=GL_m\times GL_n\subset GL_{m+n}$.  Denote $L_m:=GL_m\times \{e\}, 
L'_n:=\{e\}\times GL_n$, so that $L_{m,n}=L_m\cdot 
L'_n\cong L_m\times L'_n$. 
Let $w_{m,n}:=w_{L_{m,n},G}$, and $w_{n,m}:=w_{G,L_{m,n}}=w_{m,n}^{-1}$

It is easy to see that 
\begin{equation}
\label{eq:wmn}
w_{m,n}=(s_ms_{m-1}\cdots s_1)(s_{m+1}s_m\cdots s_2)\cdots 
(s_{m+n-1}s_{m+n-2}\cdots s_n) \ .
\end{equation}
In particular, $\dim~U^{w_{n,m}}=mn$. Denote
$C_m:=s_1\cdots s_m$ and  $C'_n:=s_{m+n-1}\cdots s_m$.
Note that $({\bf U}^{w_{n,m}},\eta^{w_{n,m}})$, 
$({\bf U}^{C_m},\eta^{C_m})$  and 
$({\bf U}^{C'_n},\eta^{C'_n})$  
are unipotent $GL_{m+n}$-crystals.
  
\begin{corollary}
$~~~$

\noindent (a) The restriction 
$({\bf U}^{w_{n,m}},\eta^{w_{n,m}})|_{L_m}$ is isomorphic in 
$U_{L_m}-{\mathcal Cryst}$ to 
\begin{equation}
\prod_{l=1}^n ({\bf U}^{C_m},\eta^{C_m})|_{L_m} \ .
\end{equation}
 
\noindent (b) The restriction 
$({\bf U}^{w_{n,m}},\eta^{w_{n,m}})|_{L'_n}$  is isomorphic in 
$U_{L'_n}-{\mathcal Cryst}$ to 
\begin{equation}
\prod_{k=1}^m ({\bf U}^{C'_n},\eta^{C'_n})|_{L'_n} \ .
\end{equation}
\end{corollary}

\subsection{$W$-invariant functions on standard unipotent $G$-crystals}

\label{subsect:restriction standard to L}

For each  $w\in W$ let $\widehat U(w)$ be the set of all
$\chi\in \widehat U$ such that $\chi(\overline w^{\,-1}u\overline w)=\chi(u)$
for all $u\in U(w)$. 

Clearly, $\widehat U(e)=\widehat U(w_0)=\widehat U$.

For any $\chi\in \widehat U$ define a regular 
function $\chi^w:BwB\to \GG_a$  
by  
\begin{equation}
\label{eq:chiw}
\chi^w(g)= \overline \chi(\overline w^{\,-1}g) \ .
\end{equation}

Let $\rho^\vee\in X_\star(T^{ad})$ be the co-character of $T^{ad}=T/Z(G)$ 
such that $\left<\rho^\vee,\alpha_i\right>=1$ for $i\in I$. Define the anti-automorphism $\iota$ of $G$  by 
$\iota(g)=Ad~\rho^\vee(-1)~(g^{-1})$. 

By definition, $\iota(t)=t^{-1}$ for $t\in T$  
and $\iota \circ x_i =x_i, ~\iota \circ y_i=y_i$
for all $i\in I$.

It is easy to see that $\iota(\overline w)=\overline {w^{-1}}$ for any $w\in W$ and 
$\chi\circ \iota=\chi$ for any $\chi \in \widehat U$. 

For any $\chi,\chi'\in \widehat U$ and any $w\in W$ define a regular function  
$f_{\chi,\chi'}^w$ on $BwB$ by
\begin{equation}
\label{eq:fchi}
f_{\chi,\chi'}^w=\chi^w+{\chi'}^{w^{-1}}\circ \iota \ .
\end{equation}

Let $L=L_J$ be any standard Levi subgroup of $G$. For any 
$\chi=\sum_{i\in I} a_i\chi_i\in \widehat U$ let $\chi^L\in \widehat U$ be defined by 
$$\chi^L:=\sum_{i\in I\setminus J} a_i\chi_i \ .$$
By definition, $\chi^L|_{U_L}\equiv 0$.

\begin{theorem} 
\label{th:chi sigma L invariant}
For any $\chi\in \widehat U(w_{L,G})$ and any $t\in Z(L)$  
the restriction of the function $f_{\chi,\chi^L}^{w_{G,L}}$ to $t\cdot B^-_{w_{L,G}}$ 
is invariant under the unipotent action of $W$ on $t\cdot B^-_{w_{L,G}}$.
\end{theorem}

Next, we describe for each $w\in W$ a basis for the subspace 
$\widehat U(w)\subset \widehat U$. 
For each $w\in W$ define a bijective map 
$\zeta_w:I(w)\to I(w^{-1})$ by 
$w(\alpha_i)=\alpha_{\zeta_w(i)}$.  
This  $\zeta_w$ is a partial bijection $I\to I$ 
(see the Appendix) with $dom(\zeta)=I(w)$ and $ran(\zeta)=I(w^{-1})$. 

For each $i\in I$ we define the set ${\bf o}_w(i)$ as follows. 
For $i\in dom(\zeta) \cup ran(\zeta)$ we put
$${\bf o}_w(i):=
\{\ldots,(\zeta_w)^{-2}(i),(\zeta_w)^{-1}(i),i,\zeta_w(i),
(\zeta_w)^2(i),\ldots \}$$ 
(where for each $k\in \ZZ$ the power $(\zeta_w)^k$ is considered to be a partial 
bijection $I\to I$).
For $i\notin dom(\zeta) \cup ran(\zeta)$ we define ${\bf o}_w(i):=\{i\}$.

For any $J\subset I$ let $\chi_J\in \widehat U$ be defined by $\chi_{J}=
\sum_{j\in J}\chi_j$.

\begin{lemma} For each $i\in I$ the function $\chi_{{\bf o}_w(i)}$ 
belongs to $\widehat U(w)$, and the set of all 
$\chi_{{\bf o}_w(i)}$, $i\in I$, is a 
basis in the vector space $\widehat U(w)$. 

\end{lemma}

Recall that the unipotent action of $W_J$ was defined in 
section \ref{subsect:adjoint model}.
 
\begin{theorem}
\label{th:W invariants}
$~~~$

\noindent (a) For each $\chi\in \widehat U(w_{G,L})$ the restriction of 
$\chi$ to $U^{w_{G,L}}$ is 
invariant under the unipotent action of $W_J$ on the 
unipotent $L$-crystal 
$({\bf U}^{w_{G,L}}, \eta^{w_{G,L}})|_L$. 

\noindent (b) For each $\chi\in \widehat U(w_{L,G})$ the restriction of 
the function  $\chi^{w_{L,G}}$ to $B^-_{w_{L,G}}$ is 
invariant under the unipotent action of $W_J$ on the 
unipotent $L$-crystal $({\bf B}^-_{w_{L,G}}, {\rm id}_{w_{L,G}})|_L$.   
\end{theorem}

\subsection{Positive structures on standard unipotent crystals}

\label{subsect:positive unipotent structures}

For any $i\in I$ define the morphism 
$\pi_i:\GG_m\to B^-$ by the formula
$\pi_i(c):= y_i\left(\frac{1}{c}\right)\alpha_i^\vee(c)$. 
Clearly, $\pi_i$ is a biregular isomorphism 
$\GG_m\widetilde \to B^-_{s_i}$.   The isomorphism $\pi_i:\GG_m\widetilde \to B^-_{w_i}$ 
defines a structure of 
$U$-variety on $\GG_m$ and a unipotent crystal on 
$\GG_m$ which we denote by $(\GG_m,\pi_i)$. 

For each sequence $\ii=(i_1,\ldots,i_l)$  
define the regular morphism
$\pi_\ii:(\GG_m)^l\to B^-$ by the formula 
\begin{equation}
\label{eq:pi general}
\pi_\ii(c_1,\ldots,c_l)=\pi_{i_1}(c_1)\cdots \pi_{i_l}(c_l) \ .
\end{equation}

For any sequence $\ii=(i_1,\ldots,i_l)\in I^l$ define 
$w_{\star}(\ii):=s_{i_1}{\star}\cdots {\star} s_{i_l}$. 

\begin{proposition}

\label{pr:factorization B-}
 For any sequence 
$\ii=(i_1,\ldots,i_l)\in I^l$ we have:

\noindent (a) There exists a sub-torus $\widetilde T_\ii\subset T$ such that  
$[\pi_\ii]$ is a dominant rational morphism  
$(\GG_m)^l\to B^-_{w_{\star}(\ii)}\cdot \widetilde T_\ii$.
The morphism $[\pi_\ii]$ is a birational isomorphism if and only if 
\begin{equation}
\label{eq:isomorphism piii}
l=l(w_{\star}(\ii))+\dim  \widetilde T_\ii \ .
\end{equation}

\noindent (b) The rational morphism 
$[\pi_\ii]$ induces the morphism in 
$U-{\mathcal Cryst}$: 
$$(\GG_m,\pi_{i_1})\times \cdots \times 
(\GG_m,\pi_{i_l})\to 
\left({\bf B}^-_{w_{\star}(\ii)}\cdot \widetilde T_\ii,
{\rm id}_{B^-_{w_{\star}(\ii)}\cdot \widetilde T_\ii}\right) \ .$$

\noindent (c) If the equality (\ref{eq:isomorphism piii}) holds 
then $[\pi_\ii]$
is a positive structure on the 
unipotent crystal 
$\left({\bf B}^-_{w_{\star}(\ii)}\cdot \widetilde T_\ii,
{\rm id}_{B^-_{w_{\star}(\ii)}\cdot \widetilde T_\ii}\right)$.

\noindent (d) There is an inductive formula for the computation 
of $\widetilde T_\ii$:

$\widetilde T_{(i)}=\alpha_i^\vee(\GG_m)$, and 
for any sequence $\ii=(i_1,\ldots,i_l)\in I^l$ and any $i\in I$ one has 
$$\widetilde T_{(\ii,i)}=
\begin{cases}
s_i(\widetilde T_\ii), & \text{if $l(w_{\star}(\ii)s_i)=l(w_{\star}(\ii))+1$,} \\
s_i(\widetilde T_{\ii})\cdot \widetilde T_i, & \text{if $l(w_{\star}(\ii)s_i)=
l(w_{\star}(\ii))-1$.} 
\end{cases} 
$$
\end{proposition}

\noindent {\bf Remarks}. 

\noindent 1. For each $w\in W$, $\ii\in R(w)$, 
the tropicalization of the  geometric crystal 
induced by  
$({\bf B}^-_w,{\rm id}_w)$ with respect to $\theta'_\ii$ is equal to the 
corresponding free combinatorial crystal ${\mathcal B}_{\ii}$ which was
constructed in \cite{k93} as a 
product of $1$-dimensional crystals 
${\mathcal B}_{i_1},\ldots, {\mathcal B}_{i_l}$. 
By the construction, ${\mathcal B}_\ii$ is a free 
$W_{|\ii|}$-crystal, where $|\ii|=\{i_1,\ldots,i_l\}$.  

\noindent 2. For any dominant
$\lambda^\vee\in \Lambda^\vee$ the image of Kashiwara's embedding 
(\ref{eq:Kashiwara embedding}) can be described as follows (see e.g. \cite{bz3}).  
Let $\chi=\sum_{i\in I} \chi_i\in \widehat U$. For any $\ii\in R(w_0)$ let 
$\theta_\ii=id_T\times [\pi_\ii]: T\times (\GG_m)^{l_0}\widetilde \longrightarrow T\cdot 
B^-_{w_0}$
be the positive structure on $T\cdot B^-_{w_0}$. Then the image of the embedding
$B(V_{\lambda^\vee})\hookrightarrow B_\ii=\ZZ^{l_0}$ is
$$\{b\in B_\ii:~Trop_{\theta_\ii}\left(f_{\chi,\chi}^{w_0}\right)(\lambda^\vee,b)\ge 0\} \ .$$

\smallskip

For each $\ii=(i_1,\ldots,i_l)\in I^l$  
define the regular morphism
$\pi^\ii:(\GG_m)^l\to U$ by 
\begin{equation}
\label{eq:pi unipotent general}
\pi^\ii(c_1,\ldots,c_l)=x_{i_1}(c_1)\cdots x_{i_l}(c_l) \ .
\end{equation}

\begin{proposition} 
\label{pr:1-dim product unipotent bruhat}
For any sequence $\ii=(i_1,\ldots,i_m)\in I^l$ we have:

\noindent (a)  The morphism $\pi^\ii$ induces a dominant rational 
morphism 
$(\GG_m)^l\to U^{(w_{\star}(\ii))^{-1}}$.  
If $\ii\in R(w)$ for some $w\in W$ then $\pi^\ii$ is 
an open inclusion $(\GG_m)^{l(w)}\hookrightarrow U^w$.

\noindent (b) If the sequence  $\ii$ belongs to $R(w)$ for some $w\in W$ 
then the birational isomorphism 
$[\pi^\ii]:(\GG_m)^{l(w)}\widetilde \to U^{w^{-1}}$ is a positive structure 
on the unipotent crystal $({\bf U}^{w^{-1}},\eta^{w^{-1}})$.
\end{proposition}

\subsection{Duality and symmetries for standard unipotent crystals}

It is easy to see that the inverse
${\cdot}^{-1}:B^-\to B^-$ induces the isomorphism 
$B^-_w\widetilde \to B^-_{w^{-1}}$
for each $w\in W$.  

\begin{lemma} The inverse in $B^-$ induces the 
isomorphism of unipotent crystals
$$({\bf B}^-_w,{\rm id}_w)^*\widetilde \to
({\bf B}^-_{w^{-1}},{\rm id}_{w^{-1}})$$
for each $w\in W$. In particular, the unipotent crystal 
$({\bf B}^-_w,{\rm id}_w)$ is self-dual if $w^2=e$. 
\end{lemma}

\section{Proofs}

\label{sect:proofs}

\subsection {Proofs of results in section \ref{sect:definitions}}
$ ~~~$

\noindent {\bf Proof of Theorem  \ref{th:X(infinity)}}.
Since theorem \ref{th:X(infinity)} is a particular case of Theorem 
\ref{th:the first functor} we refer to the proof of 
Theorem \ref{th:the first functor} in 
section \ref{subsect:proofs section 3}.
\hfill$\square$\smallskip

\noindent {\bf Proof of Proposition \ref{pr:W-action}}.
Clearly, $s_i:X\to X$ is a birational involution. 
Thus, it suffices to prove is that for any $w_1,w_2\in W$ satisfying
$l(w_1w_2)=l(w_1)+l(w_2)$ one has:  
$$w_1(w_2(x))=(w_1w_2)(x) \ .$$

As it follows from definition of 
$e_w:T\times X\to X$ we have
$$e_{w_1w_2}^t=e_{w_1}^{w_2(t)}\circ e_{w_2}^t$$
for all $w_1,w_2\in W$ satisfying $l(w_1w_2)=l(w_1)+l(w_2)$.  
Note that for any $w\in W$ we have 
$$\gamma(e_w^t(x))=t\cdot w(t^{-1})\cdot \gamma(x) \ .$$
This implies that $\gamma(w(x))=\gamma(e_w^{\gamma(x)^{-1}}(x))=w(\gamma(x))$. 
Therefore, for any $w_1,w_2\in W$ such that  $l(w_1w_2)=l(w_1)+l(w_2)$
we have 
$$w_1(w_2(x))=e_{w_1}^{\gamma(w_2(x))^{-1}}\left(
e_{w_2}^{\gamma(x)^{-1}}(x)\right)$$
$$=e_{w_1}^{w_2(\gamma(x)^{-1})}e_{w_2}(x)^{\gamma(x)^{-1}}=
e_{w_1w_2}^{\gamma(x)^{-1}}(x)=(w_1w_2)(x) \ .$$

Proposition \ref{pr:W-action} is proved. $\quad\square$

\noindent {\bf Proof of Theorem \ref{th:positive structure}}. 
Let $W\times T'\to T'$ be the action of $W$ 
obtained from the action \ref{eq:W-action} by twisting with $\theta$, 
that is, $(w,t')\mapsto \theta^{-1}(w(\theta(t'))$. It is easy to see that 
$W$ acts on $T'$ by a positive birational isomorphisms. 
Therefore, applying the functor $Trop$ to the action $W\times T'\to T'$ 
we obtain an action of $W$ on $X_\star(T')$. 

Theorem \ref{th:positive structure} is proved. 
\hfill$\square$\smallskip

\noindent {\bf Proof of Theorem \ref{th:trivialization tau}}. 
By definition, for $i,j\in I=\{1,2,\ldots,r\}$ we have: 
$$\omega_i(e_j^c(x))=c^{\delta_{ij}}\omega_i(x) \ .$$

Thus $\tau(e_1^c(x))=\tau(x)$ and therefore $\tau(s_1(x))=\tau(x)$. 

Let us compute $\tau(e_j^c(x))$ for $j>1$. The computation is based on
the following  obvious statement. 
\begin{lemma} For each $j=2,\ldots,r$ we have:
\begin{equation}
\label{eq:split tau}
\tau(e_j^c(x))=(e_{(j;r)}^x\circ \tau_j)\left(e_j^c(x)\right) \ ,
\end{equation}
where 
\begin{equation}
\label{eq:tauj}
\tau_j(z)=\left(e_1^{\frac{\omega_{j+1}(z)}{\omega_j(z)}}\cdots 
e_j^{\frac{\omega_{j+1}(z)}{\omega_j(z)}}\right)
\cdots \left(e_1^{\frac{\omega_3(z)}{\omega_{2}(z)}}
e_2^{\frac{\omega_3(z)}{\omega_{2}(z)}}\right)
\left(e_1^{\frac{\omega_2(z)}{\omega_1(z)}}\right)(z) \ ,
\end{equation}
$$e_{(j;r)}^x=\left(e_1^{\frac{\omega_{r+1}(x)}{\omega_r(x)}}\cdots 
e_r^{\frac{\omega_{r+1}(x)}{\omega_r(x)}}\right)\cdots
\left(e_1^{\frac{\omega_{j+2}(x)}{\omega_{j+1}(x)}}\cdots 
e_{j+1}^{\frac{\omega_{j+2}(x)}{\omega_{j+1}(x)}}\right)$$
and we use the convention $\omega_{r+1}(x)\equiv 1$ and 
$e_{(r;r)}^x\equiv 1$. 

\end{lemma}

Substituting $c=\frac{1}{\alpha_j(\gamma(x))}
=\frac{\omega_{j-1}(x)\omega_{j+1}(x)}
{\omega_j(x)}$ into (\ref{eq:split tau}) 
we obtain 
$$\tau(s_j(x))=e_{(j;r)}(x)\circ \tau_j(s_j(x)) \ .$$ 
This implies that all we have to do to prove  
Theorem \ref{th:trivialization tau} is to check that 
\begin{equation}
\label{eq:desired identity}
\tau_j(s_j(x))=\tau_j(x)
\end{equation} 
for $j=2,\ldots,r$. 

Let us compute first  $\tau_j(e_j^c(x))$, $j\ge 2$.  Substituting  
$z=e_j^c(x)$ in (\ref{eq:tauj}) and using the equalities 
$\omega_k(z)=c^{\delta_{jk}}\omega_j(x)$ and  $\tau_{j-2}\circ e_j^c
=e_j^c \circ \tau_{j-2}$ we obtain
$$\tau_j(e_j^c(x))=
\left(e_1^{\frac{\omega_{j+1}(x)}{c\omega_j(x)}}\cdots 
e_j^{\frac{\omega_{j+1}(x)}{c\omega_j(x)}}\right)
\left(e_1^{\frac{c\omega_j(x)}{\omega_{j-1}(x)}}\cdots 
e_{j-1}^{\frac{c\omega_j(x)}{\omega_{j-1}(x)}}\right)
\left(e_j^c\circ \tau_{j-2}(x)\right) \ . $$

Substituting $c=\frac{\omega_{j-1}(x)\omega_{j+1}(x)}{\omega_j(x)}$ into
(\ref{eq:tauj}), we obtain 
$$\tau_j(s_j(x))=\left(e_1^{\frac{\omega_j(x)}{\omega_{j-1}(x)}}\cdots 
e_j^{\frac{\omega_j(x)}{\omega_{j-1}(x)}}\right)
\left(e_1^{\frac{\omega_{j+1}(x)}{\omega_j(x)}}\cdots 
e_{j-1}^{\frac{\omega_{j+1}(x)}{\omega_j(x)}}\right)
e_j^{\frac{\omega_{j-1}(x)\omega_{j+1}(x)}
{\omega_j(x)}}( \tau_{j-2}(x)) \ .$$

In order to finish the proof of Theorem \ref{th:trivialization tau} 
we will use the following result.

\begin{lemma} 
\label{le:the relation}
For any $j=1,2,\ldots,r$ the following relation holds.
\begin{equation} 
\label{eq:crystal relation}
(e_1^c\cdots e_j^c)(e_1^{c'}\cdots e_{j-1}^{c'})e_j^{\frac{c'}{c}}=
(e_1^{c'}\cdots e_j^{c'})(e_1^c\cdots e_{j-1}^c)
\end{equation}

\end{lemma}

\begin{proof} 
Induction in $j$. For $j=1$ the identity
$e_1^c e_1^{\frac{c'}{c}}=e_1^{c'}$ is true.  

Now let $j>1$. Using the commutation of $e_j$ with each of 
$e_1,\ldots,e_{j-2}$ we rewrite the left hand side of 
\ref{eq:crystal relation} as:
$$(e_1^c\cdots e_{j-1}^c)(e_1^{c'}\cdots 
e_{j-2}^{c'})e_j^ce_{j-1}^{c'}e_j^{\frac{c'}{c}} \ .$$ 
Applying the to above expression the basic relation (\ref{eq:A2}) 
written in the form 
$e_j^ce_{j-1}^{c'}e_j^{\frac{c'}{c}}=
e_{j-1}^{\frac{c'}{c}}e_j^{c'}e_{j-1}^c$, 
we see that the left hand side of 
\ref{eq:crystal relation} equals:
$$(e_1^{c'}\cdots e_j^{c'})(e_1^c\cdots 
e_{j-1}^c)e_{j-1}^{\frac{c'}{c}}e_j^{c'}e_{j-1}^c$$ 
Finally, applying the inductive hypothesis 
(\ref{eq:crystal relation}) with $j-1$, we obtain 
$$(e_1^{c'}\cdots e_{j-1}^{c'})(e_1^c\cdots e_{j-2}^c)
e_j^{c'}e_{j-1}^c=(e_1^{c'}\cdots e_j^{c'})
(e_1^c\cdots e_{j-1}^c) \ .$$ 

Lemma is proved.   
\end{proof}

We see that  (\ref{eq:desired identity}) is a special case of 
Lemma \ref{le:the relation} with $c=\frac{\omega_j(x)}{\omega_{j-1}(x)},
c'= \frac{\omega_{j+1}(x)}{\omega_j(x)}$. 
 
Theorem \ref{th:trivialization tau} is proved. 
\hfill$\square$

\subsection{Proof of results in section \ref{sect:unipotent crystals}}

\label{subsect:proofs section 3}

$~~~~$

\noindent {\bf Proof of Theorem \ref{th:product}}. Prove (a). 
We start with the following result.

\begin{lemma} For any unipotent crystal $({\bf X},{\bf f})$ 
the following identity holds. 
\begin{equation}
\label{eq:cocycle}
\pi(u'\cdot {\bf f}(u(x)))\cdot \pi(u\cdot {\bf f}(x))=
\pi((u'u)\cdot {\bf f}(x))
\end{equation}
for any $u,u'\in U$, $x\in X$. 
\end{lemma}

\begin{proof} Denote $b={\bf f}(x)$. Note that 
${\bf f}(u(x))=u(b)=u\cdot b\cdot \pi(u\cdot b)^{-1}$. 
Then the identity  (\ref{eq:cocycle}) can be rewritten as
$$\pi(u'\cdot u\cdot b\cdot \pi(u\cdot b)^{-1})\cdot \pi(u\cdot b)=
\pi(u'u\cdot b)$$
This identity is true because the morphism $\pi:G\to U$ is right 
$U$-equivariant. Lemma is proved. \end{proof}  

One can easily see that (\ref{eq:cocycle}) implies that 
$(u'u)(x,y)=u'(u(x,y))$ for all  $u,u'\in U$, $(x,y)\in X\times Y$. 

Thus, the formula 
(\ref{eq:unipotent group action on the product}) defines an $U$-action. 

Part (a)  is proved. Prove (b). We need the following fact.

\begin{lemma} 
\label{le:the property of pi}
For $u\in U$, 
$b,b'\in B^-$ we have
\begin{equation}
\label{eq:the property of pi}
\pi\left(u\cdot b\cdot b'\right)=
\pi\left(\pi(u\cdot b)\cdot b'\right)
\end{equation}
\end{lemma}

\begin{proof} 
It suffices to take $u=x_i(a)$. In this case Lemma 
\ref{le:elementary pi} implies that 
for $u,u'\in U^-, t,t'\in T$ we have:
$$\pi(x_i(a)\cdot  u\cdot tu'\cdot t')=
\pi(x_i(a)\cdot  u\cdot tu't^{-1}\cdot tt')
=x_i\left((a^{-1}+\chi^-_i(utu't^{-1}))^{-1}\cdot 
\alpha_i(tt')^{-1})\right)$$
$$=x_i\left((a^{-1}+\chi^-_i(u)+\alpha_i(t^{-1})\chi_i(u'))^{-1}\cdot 
\alpha_i(t^{-1})\cdot \alpha_i(t^{'-1})\right)$$
$$=x_i\left({a'}^{-1}+\chi^-_i(u'))^{-1}\alpha_i(t^{'-1})\right)=
\pi(x_i(a')\cdot u't')$$
where $a'=(a^{-1}+\chi^-_i(u))^{-1}\cdot \alpha_i(t^{-1})$. 
Lemma is proved. \end{proof}

Part (b) is proved, and Theorem \ref{th:product} is also proved. 
\hfill$\square$\smallskip

\noindent {\bf Proof of Proposition \ref{pr:monoidal category}}. 
Let $({\bf X}_1,{\bf f}_{X_1})$, 
$({\bf X}_2,{\bf f}_{X_2}), ({\bf X}_3,{\bf f}_{X_3})$ 
be unipotent crystals and let 
$$({\bf X}_{12,3},{\bf f}_{X_{12,3}}):=\left(({\bf X}_1,{\bf f}_{X_1})
\times ({\bf X}_2,{\bf f}_{X_2})\right)\times  
({\bf X}_3,{\bf f}_{X_3}) \ ,$$
$$ ({\bf X}_{1,23},{\bf f}_{X_{1,23}}):=({\bf X}_1,{\bf f}_{X_1})
\times \left(({\bf X}_2,{\bf f}_{X_2})\times  
({\bf X}_3,{\bf f}_{X_3})\right) \ .$$ 
  
For any $x_k\in X_k$, $k=1,2,3$ we have
$$ {\bf f}_{X_{12,3}}((x_1,x_2),x_3)={\bf f}_{X_1}(x_1)
\cdot {\bf f}_{X_2}(x_2)\cdot 
{\bf f}_{X_3}(x_3)={\bf f}_{X_{1,23}}(x_1,(x_2,x_3)) \ ,$$
that is, ${\bf f}_{X_{12,3}}={\bf f}_{X_{1,23}}$. 

It suffices to prove that ${\bf X}_{12,3}={\bf X}_{1,23}$, that is, these 
two $U$-actions on $X_1\times X_2\times X_3$ are equal. 

For $x_k\in X_k$, $k=1,2,3$ denote $b_k:={\bf f}_{X_k}(x_k)$. 

Let us write the $U$-actions respectively on
${\bf X}_{12,3}$, ${\bf X}_{1,23}$: 
$$u((x_1,x_2),x_3)=
(u(x_1,x_2),   \pi\left(u\cdot b_1\cdot b_2\right)(x_3)  )$$
$$=(u(x_1),\pi(u\cdot b_1)(x_2),
\pi\left(u\cdot b_1\cdot b_2\right)(x_3)) \ .$$
$$u(x_1,(x_2,x_3))=(u(x_1),\pi(u\cdot b_1)(x_2,x_3))$$
$$=(u(x_1),\pi(u\cdot b_1)(x_2),\pi\left(\pi(u\cdot b_1)\cdot
b_2\right)(x_3)) \ .$$

It follows from Lemma \ref{le:the property of pi} that 
$u((x_1,x_2),x_3)=u(x_1,(x_2,x_3))$. 

This proves Proposition \ref{pr:monoidal category}.
\hfill$\square$\smallskip

\noindent {\bf Proof of Proposition \ref{pr:G is crystal}}. 
It follows from the definition (\ref{eq:gauss}) of ${\bf g}$ and the definition 
of $\pi$ that for $g\in G, u\in U$ 
one has: 
$${\bf g}(u\cdot g)={\bf f}_G(u\cdot g))\cdot \pi(u\cdot g) \ .$$
Note that ${\bf f}_G(u\cdot g)=u\cdot \left({\bf f}_G(g)\right)$ 
and 
$\pi(ug)=\pi\left(u\left({\bf f}_G(g)\right)\right)\cdot \pi(g)$. 
Thus, Proposition \ref{pr:G is crystal} follows from Theorem
\ref{th:product}.  
\hfill$\square$\smallskip

\noindent {\bf Proof of Theorem \ref{th:the first functor}}
In the view of Lemma \ref{le:verma}, it suffices to 
prove the relations (\ref{eq:A1A1})-(\ref{eq:G2}).

First, we are going to prove the relations (\ref{eq:A1A1}) 
and (\ref{eq:A2}). For $\varepsilon\in \{0,-1\}$ let 
$$(I\times I)_\varepsilon=
\{(i,j)\in I\times I: i\ne j~{\rm and }~
\left<\alpha_i^\vee,\alpha_j\right>= 
\left<\alpha_j^\vee,\alpha_i\right>= \varepsilon\} \ .$$
Note that if $G$ is simply-laced, then 
$(I\times I)_0 \cup (I\times I)_{-1}\cup \Delta(I)=I\times I$.

\begin{lemma} 
$~~~$

\noindent (a) For any  $(i,j)\in (I\times I)_0$ we have
\begin{equation}
\label{eq:commutativity}
x_i(a_1)x_j(a_2)=x_j(a_2)x_i(a_1)  \ .
\end{equation}
\noindent (b) For any $(i,j)\in (I\times I)_{-1}$ we have 
\begin{equation}
\label{eq:braid relation1}
x_i(a_1)x_j(a_2)x_i(a_3)=x_j\left(\frac{a_2a_3}{a_1+a_3}\right)x_i
\left(a_1+a_3\right)x_j\left(\frac{a_1a_2}{a_1+a_3}\right)
\end{equation} 
\end{lemma}

\begin{proof} Clear. \end{proof}

\begin{proposition} Let $({\bf X},{\bf f})$ be a unipotent 
$G$-crystal and let ${\mathcal X}_{ind}$ be the induced geometric crystal. 
For any $i,j\in (I\times I)_0 \cup (I\times I)_{-1}$ such that 
$\varphi_i\not \equiv 0,\varphi_j\not \equiv 0$, one has:

\noindent (a) if  $(i,j)\in (I\times I)_0 $ then 
\begin{equation} 
\varphi_i(e_j^c(x))=\varphi_i(x), \varphi_i(e_j^c(x))=\varphi_i(x) \ ;
\end{equation}
\noindent (b) If  $(i,j)\in (I\times I)_{-1}$ 
then there exist rational functions 
$\varphi_{ij}, \varphi_{ij}$ on $X$ such that:
\begin{equation}
\label{eq:canonical basis}
 \varphi_i(x)\varphi_j(x)=\varphi_{ij}(x)+\varphi_{ji}(x)
\end{equation}
and
\begin{equation}
\label{eq:invariancy}
\begin{matrix}
\varphi_{ij}(e_j^c(x))=\varphi_{ij}(x),
~\varphi_{ij}(e_i^c(x))=c^{-1}\varphi_{ij}(x), \\ 
\varphi_{ji}(e_i^c(x))=\varphi_{ji}(x),  
~\varphi_{ji}(e_j^c(x))=c^{-1}\varphi_{ji}(x)
\end{matrix}
\end{equation}

\noindent (c) We have
\begin{equation}
\label{eq:ansatz relation1}
 \varphi_j(x)\varphi_i(e_j^c(x))=c\varphi_{ij}(x)+\varphi_{ji}(x), 
~~\varphi_i(x)\varphi_j(e_i^c(x))=c\varphi_{ji}(x)+\varphi_{ij}(x)
\end{equation}
\begin{equation}
\label{eq:ansatz relation2}
\begin{matrix}
\varphi_i(e_j^{c_1}e_i^{c_2}(x))=
\varphi_i(x)\frac{c_1c_2^{-1}\varphi_{ij}(x)+\varphi_{ji}(x)}
{c_2\varphi_{ji}(x)+\varphi_{ij}(x)},
\varphi_j(e_i^{c_1}e_j^{c_2}(x))=
\varphi_j(x)\frac{c_1c_2^{-1}\varphi_{ji}(x)+\varphi_{ij}(x)}
{c_2\varphi_{ij}(x)+\varphi_{ji}(x)}
\end{matrix}
\end{equation}

\end{proposition}

\begin{proof} It suffices to prove the proposition in the assumption that 
$I=\{i,j\}$, that is, when $G$ is semisimple of types $A_1\times A_1$ and 
$A_2$ respectively.  The part (a) follows. 

Let us prove (b). It suffices to analyze the case when  $G=GL_3$. Due to Lemma 
\ref{le:universal G}  it suffices to prove the statement  only for  $X=G$, ${\bf f}_X=\pi^-$.  
In this case $I=\{i,j\}$, and we set $i:=1, j:=2$ in the standard way.

It is easy to see that 
$\varphi_k(g)=\overline \chi_k^{\,-}(g)=\frac{\Delta'_k(g)}{\Delta_k(g)}$ for $k=1,2$, where
$$\Delta_1(g)=g_{11}, \Delta'_1(g)=g_{21},
~\Delta_2(g)=\det\begin{pmatrix}  g_{11} &  g_{12} \\
              g_{21} & g_{22}    \\
      \end{pmatrix},~~ 
\Delta'_2(g)=\det\begin{pmatrix}  g_{11} &  g_{12} \\
              g_{31} & g_{32}    \\
      \end{pmatrix}\ .$$
Furthermore,  is easy to see that the actions 
$e_1,e_2:G_m\times B^-\to B^-$ are given by the formula 
(\ref{eq:simple multiplicative generator ei on G}):
$$e_k^c(g)=x_k\left((c-1)\frac{\Delta_k(g)}{\Delta'_k(g)}\right)\cdot g$$
for $k=1,2$. 

Define the functions 
$\varphi_{12},\varphi_{21}$ on $G$ by:
$$\varphi_{12}(g)=\frac{\Delta''_1}{\Delta_1}, 
\varphi_{21}(g)=\frac{\Delta''_2(g)}{\Delta_2(g)} \ ,$$
where
$$\Delta''_1(g)=g_{31},
~\Delta''_2(g)=\det\begin{pmatrix}  g_{21} &  g_{22} \\
              g_{31} & g_{32}    \\
      \end{pmatrix}$$

It is easy to see that 
$$\Delta'_1(g)\Delta'_2(g)=\Delta_1(g)\Delta''_2(g)+\Delta''_1(g)\Delta_2(g) \ .$$
This identity implies (\ref{eq:canonical basis}).  Furthermore, 
it is easy to see that 
$$\Delta_k(e_l^c(g))=c^{\delta_{kl}}\Delta_k(g),~\Delta''_k(e_l^c(g))=\Delta''_k(g)$$
for $k,l\in \{1,2\}$ and $\Delta'_k(e_k^c(g))=\Delta_k(g)$ for $k=1,2$. 
This implies (\ref{eq:invariancy}). Part (b) is proved.

Part (c) easily follows from (b).  \end{proof}

In order to  prove the relations (\ref{eq:A1A1}) for any 
$(i,j)\in (I\times I)_0$ we compute the left hand side and the right hand 
side of (\ref{eq:A1A1}). By definition
$$e_i^{c_1}e_j^{c_2}(x)=
x_i(a_1)x_j(a_2)(x) \ ,$$
where $a_1=\frac{c_1-1}{\varphi_i(e_j^{c_2}(x))}
=\frac{c_1-1}{\varphi_i(x)}$, 
and $a_2=\frac{c_2-1}{\varphi_j(x)}$. Analogously
$$e_j^{c_2}e_i^{c_1}(x)=x_j(a_2)x_i(a_1)(x) \ .$$
Thus, using the relation (\ref{eq:commutativity}), we see that 
the left and the right hand sides 
of (\ref{eq:A1A1}) are equal. This proves all the relations (\ref{eq:A1A1}). 
 
To prove the relations (\ref{eq:A2}) for any 
$(i,j)\in (I\times I)_{-1}$ we compute the left hand side and the right
hand side of (\ref{eq:A1A1}).  By definition
$$e_i^{c_1}e_j^{c_1c_2}e_i^{c_2}(x)=x_i(a_1)x_j(a_2)x_i(a_3)(x) \ , $$
where 
$$a_3=\frac{c_2-1}{\varphi_i(x)},a_2=
\frac{c_1c_2-1}{\varphi_j(e_i^{c_2}(x))}
=\frac{\varphi_i(x)(c_1c_2-1)}{c_2\varphi_{ji}(x)+\varphi_{ij}(x)}, $$
and 
$$a_1=\frac{c_1-1}{\varphi_i(e_j^{c_1c_2}e_i^{c_2}(x))}=
\frac{(c_1-1)(c_2\varphi_{ji}(x)+\varphi_{ij}(x))} 
{\varphi_i(x)(c_1\varphi_{ij}(x)+\varphi_{ji}(x))} \ .  $$
Similarly, 
$$e_j^{c_2}e_i^{c_1c_2}e_j^{c_1}(x)=x_j(a'_1)x_i(a'_2)x_j(a'_3)(x)  \ ,$$
where 
$$a'_3=\frac{c_1-1}{\varphi_j(x)},a'_2=
\frac{c_1c_2-1}{\varphi_i(e_j^{c_1}(x))}
=\frac{\varphi_j(x)(c_1c_2-1)}{c_1\varphi_{ij}(x)+\varphi_{ji}(x)}, $$
and 
$$a'_1=\frac{c_2-1}{\varphi_j(e_i^{c_1c_2}e_j^{c_1}(x))}=
\frac{(c_2-1)(c_1\varphi_{ij}(x)+\varphi_{ji}(x))} 
{\varphi_j(x)(c_2\varphi_{ji}(x)+\varphi_{ij}(x))} \ .  $$

It is easy to see that
$$a'_1a'_2=a_2a_3, a'_2a'_3=a_1a_2, a'_2=a_1+a_3 \ .$$
Thus it follows from (\ref{eq:braid relation1}) that
$$e_i^{c_1}e_j^{c_1c_2}e_i^{c_2}(x)=x_i(a_1)x_j(a_2)x_i(a_3)(x)
=x_j(a'_1)x_i(a'_2)x_j(a'_3)(x)=e_i^{c_1}e_j^{c_1c_2}e_i^{c_2}(x)$$

This proves all the relations (\ref{eq:A2}). 

\smallskip

Thus, we have proved Theorem \ref{th:the first functor} for all 
simply-laced reductive groups $G$. In particular, 
the relations (\ref{eq:verma}) hold for such groups. 

It remains to prove the relations (\ref{eq:B2}) and (\ref{eq:G2}). 
Instead of doing the computations directly, we  will prove the relations
(\ref{eq:verma}) for any group $G$ by deducing  these relations from  
the relations (\ref{eq:verma}) for a certain simply-laced group 
$G'$ containing $G$. 

Without loss of generality we may assume that $G$ is adjoint semisimple.  

It is well-known that there is an  adjoint semisimple simply-laced group $G'$, 
an outer automorphism $\sigma:G'\to G'$ and an injective group homomorphism 
$f:G\hookrightarrow  G'$ such $f(G)=(G')^\sigma$.  
In what follows we identify $G$ with it's image $f(G)$.   

Moreover, one can always choose an
outer automorphism $\sigma$ in such a way that  $\sigma$ preserves a chosen
Borel subgroup $B'\subset  G'$ and a maximal torus $T'\subset B'$,  and
satisfies $B=(B')^\sigma, T=(T')^\sigma$. Let  $(B')^-$ be the 
opposite Borel subgroup containing  $T'$. Then  $(B')^-$ is also
$\sigma$-invariant and $B^-=((B')^-)^\sigma$. Also the automorphism $\sigma$ 
induces an  injection $\sigma_{\star}:\Lambda^\vee\to (\Lambda')^\vee$.
 
Let $I'$ be the Dynkin diagram of $G'$. It is easy to see that for any 
$i\in I$ there exists a subset $\tau(i)\in I'$ such that 
$\alpha_i^\vee=\sum_{i'\in \tau(i)} {\alpha'}^{\vee}_{i'}$ 
and that for any  $i',j'\in \tau(i)$ the subgroups  $U_{i'}, U_{j'}$ 
commute. This implies that the Weyl group  $W$ is the subgroup of the 
Weyl group $W'$ of $G'$ and its generators $s_i$, $i\in I$, are given by: 
$$s_i=\prod_{i'\in \tau(i)} s'_{i'}  \ .$$

One can always choose the homomorphisms $x'_{i'}:
\GG_m\to U'_{i'}$, $y'_{i'}:\GG_m\to U'_{i'}$, 
$i'\in I'$ in such a way that 
$$x_i(a)=\prod_{i'\in \tau(i)} x'_{i'}(a), 
~y_i(a)=\prod_{i'\in \tau(i)} y'_{i'}(a)$$ 
for any $i\in I$. This implies that 
$\overline s_i=\prod_{i'\in \tau(i)} \overline {s'}_{i'}$.

Due to Lemma \ref{le:universal G} we may assume without loss of 
generality that $X=G$.  
Let ${\mathcal X}_{ind}$ be the geometric $G$-crystal induced by 
$({\bf G},{\rm id}_{G})$ 
and let ${\mathcal X}'_{ind}$ be the geometric $G'$-crystal induced by 
$({\bf G}',{\rm id}_{G'})$. We denote by 
$e_i:\GG_m\times G\to G$, $i\in I$, the actions 
(\ref{eq:simple multiplicative generator ei on G}) of $\GG_m$, and by  
$e'_{i'}:\GG_m\times G'\to G'$, $i'\in I'$   
the corresponding actions of $\GG_m$. As it follows from (\ref{eq:A1A1}), 
the transformations 
${e'_{i'}}^{c}$ commute with ${e'_{j'}}^{c'}$ for any $i\in I,i',j'\in \tau(i)$.

For each $i\in I$ let us fix a linear ordering $\overline \tau(i)$ of
$\tau(i)$. It is easy to see that for any reduced sequence
$\ii=(i_1,\ldots,i_l)\in I^l$ 
the sequence $\tau(\ii):=(\overline \tau(i_1);\ldots;\overline \tau(i_l))$
is also reduced.

Let $e_{\ii}:T\times G\to G$ and 
$e'_{\overline \tau(\ii)}:T'\times G'\to G'$ be morphisms as 
in (\ref{eq:general eii}).  

The following lemma is obvious. 
 
\begin{lemma} 
For each reduced sequence $\ii=(i_1,\ldots,i_l)\in I^l$ we have:  

\noindent (i)  $dom(e_\ii)=dom(e'_{\tau(\ii)})\cap (T\times G)$. 

\noindent (ii) $e'_{\overline \tau(\ii)}|_{T\times G}=e_{\ii}$.

\end{lemma}

The lemma implies that for any reduced sequences $\ii,\tilde \ii\in I^l$ 
satisfying $w(\ii)=w(\tilde \ii)$ we have 
$e_\ii=e'_{\overline \tau(\ii)}|_{T\times G}
=e'_{\overline \tau(\tilde \ii)}|_{T\times G}=e_{\tilde \ii}$. 

This proves all the relations  \ref{eq:verma} for geometric $G$-crystals. 

Theorem \ref{th:the first functor} is proved. 
\hfill$\square$\smallskip

\noindent{\bf Proof of Theorem \ref{th:positive unipotent}}. 
As it  follows from Lemma 
\ref{le:geometric crystal on the product}. 
that:

\noindent (a)  The morphism 
$\gamma_{X\times Y}\circ \theta_{X\times Y}:T'\times T''\to T$ 
is positive.

\noindent (b) Each  function $\varphi^{X\times Y}_i$ is positive, 

\noindent (c) Each action $e_i:\GG_m\times X\times Y\to X\times Y$ is 
positive.

Theorem \ref{th:positive unipotent} is proved. 
\hfill$\square$\smallskip

\noindent {\bf Proof of Proposition \ref{pr:dual unipotent crystal}}.  
Without loss of generality we may assume that 
$X$ is a subset of $B^-$ or even $X=B^-$ so that ${\bf f}_X$ 
is the identity map $X\to B^-$.
Then it is easy to see that $\alpha^*(u,b)=(u(b^{-1}))^{-1}$.
This implies that $\alpha^*$ is an action of $U$. Furthermore, 
${\bf f}^*(b)=b^{-1}$. Thus 
$${\bf f}^*(\alpha^*(u,b))=u(b^{-1})=u({\bf f}^*(b)) \ .$$

Proposition \ref{pr:dual unipotent crystal} is proved. 
\hfill$\square$\smallskip

\noindent {\bf Proof of Proposition \ref{pr:diagonalization}}. 
It is easy to see that the formula 
(\ref{eq:diagonalization}) 
is equivalent to  the equation
$$u{\bf v}({\bf f}(x))={\bf v}({\bf f}(u(x))\pi({\bf f}(x))$$
for $x\in B^-, u\in U$. ${\bf v}$ is 
two-sided $U$-equivariant  
this identity is equivalent to:
$${\bf v}(u\cdot{\bf f}(x))={\bf v}({\bf f}(u(x))\cdot 
\pi({\bf f}(u(x))) \ . $$
This identity is true since  
$${\bf f}(u(x))\pi({\bf f}(u(x))=
u\left({\bf f}(x)\right)\pi({\bf f}(u(x))
=u\cdot{\bf f}(x) \ .$$

Proposition \ref{pr:diagonalization} is proved. 
\hfill$\square$\smallskip

\noindent {\bf Proof of Theorem \ref{th:dual unipotent to geometric}}. 
We need the following obvious result. 

\begin{lemma} For any unipotent crystal $({\bf X},{\bf f})$ let  
$\varphi_i^*=\overline \chi_i^{\,-}(({\bf f}_X(x))^{-1})$  be the
function on $X$ defined by (\ref{eq:phi on X}) for  
the dual unipotent crystal $({\bf X},{\bf f})^*,i\in I$. 
Then $\varphi_i^*(x)=-\varphi_i(x)\alpha_i(\gamma(x))$.   
\end{lemma}
 
\smallskip

Using (\ref{eq:elementary pi}), we obtain 
$\alpha^*(x_i(a),x)=
x_i\left((a^{-1}+\varphi_i^*(x))^{-1}\alpha_i(\gamma^*(x)^{-1})\right)(x)$. 

Substituting $a=\frac{c-1}{\varphi_i^*(x)}$ 
into this expression for $\alpha^*$, we obtain 
the action of the corresponding multiplicative group of 
the induced dual crystal. On the other hand, 
simplifying the right hand side, we obtain 
$$\alpha^*\left(x_i\left(\frac{c-1}{\varphi_i^*(x)}\right),x\right)
=x_i\left(\frac{c-1}{-c\varphi_i^*(x)\alpha_i(\gamma^*(x))}\right)(x)=
x_i\left(\frac{c^{-1}-1}{\varphi_i(x)}\right)(x)=e_i^{c^{-1}}(x)$$

which proves Theorem \ref{th:dual unipotent to geometric}. \hfill$\square$\smallskip

\noindent {\bf Proof of Proposition \ref{pr:uw}}. 
Denote by $T_r$ the set of regular elements of $T$ and by  
$B^-_{rss}$ the pre-image $pr_T^{-1}$ which is the set of all regular 
semisimple elements in $B^-$. This is an open dense subset. 
Note that $B^-_{rss}$ intersects non-trivially each $B^-_{w'}t$, $t\in T_r$

We define the rational morphism ${\bf u}_w:B^-_{rss}\to U^w$ 
inductively. Define ${\bf u}_e:=pr_e$, 
$${\bf u}_{s_i}(b):=
x_i\left(\frac{1-\alpha_i(\gamma(b))}{\varphi_i(b)
\alpha_i(\gamma(b))}\right)$$
for $i\in I$, and for any $w',w''\in W$ such that $l(w'w'')=l(w')+l(w'')$, define:
\begin{equation}
\label{eq:inductive uw}
{\bf u}_{w'w''}(b):={\bf u}_{w'}({\bf u}_{w''}(b)\cdot b\cdot 
{\bf u}_{w''}(b)^{-1})\cdot {\bf u}_{w''}(b) \ .
\end{equation}

\begin{lemma} For each $w\in W$ the morphism ${\bf u}_w$ does not depend 
on the choice of the expression (\ref{eq:inductive uw})  
and satisfies (\ref{eq:unipotent w}). 
\end{lemma}

\begin{proof} It is easy to see that ${\bf u}_w$ satisfies
(\ref{eq:unipotent w}). Since the centralizer in $U$ of any
element $b\in B^-_{rss}$ is trivial we see that  ${\bf u}_w$ 
does not depend on a choice of the expression (\ref{eq:inductive uw}). 

Lemma is proved. \end{proof}

Furthermore, we can repeat the same definitions for any sub-variety of the 
form $B^-_{w'}t\cap B^-_{rss}$. In this case for any $i\in {\rm supp}(w')$ the
rational morphism  ${\bf u}_{s_i}:B^-_{w'}t\cap B^-_{rss}\to U^{s_i}$ is a
well-defined  and it is given by the analogous formula.    

Proposition \ref{pr:uw} is proved. \hfill$\square$

\subsection{Proof of results in section \ref{sect:examples}}
$~~~~~$

\noindent{\bf Proof of Proposition \ref{pr:commutative diagram bz3}}. 
Follows from  \cite{bz3}, 
Theorem 4.7 with $u=e, v=w$. 
\hfill$\square$\smallskip

\noindent {\bf Proof of Proposition \ref{pr:U-action on Uw}}. 
Let us express the action of $U$ on $U^w$ as the conjugation of the 
action of $U$ on $B^-_w$ with the isomorphism $\eta_w$:
$$u(x)=\eta_w(u(b)) \ ,$$
where $b=\eta^w(x)$. Using  (\ref{eq:the action on B-}), (\ref{eq:inverse eta}) 
and the right $U$-equivariancy of $\pi$,  we obtain: 
$$u(x)=\pi\left(\overline w\cdot \pi(u\cdot b)\cdot b^{-1}u^{-1}\right)^{-1}
=u\cdot \left(\pi(\overline w\cdot \pi(u\cdot b)\cdot b^{-1})\right)^{-1} \ . $$

Using (\ref{eq:action of U on Uw}) in the left hand side of the above identity, we obtain
$$\pi^w(u\cdot x)=\pi(\overline w\cdot \pi(u\cdot b)\cdot b^{-1})\cdot x
=\pi(\overline w\cdot \pi(u\cdot b)\cdot b^{-1}\cdot x) \ .$$

Finally, 
$$b^{-1}x=b^{-1}\eta_w(b)=b^{-1}\cdot (\pi(\overline w\cdot b^{-1}))^{-1}=
\overline w^{\,-1}\cdot \pi^-(\overline w\cdot b^{-1}) \ . $$
Proposition is proved. \hfill$\square$\smallskip

\noindent{\bf Proof of Proposition \ref{pr:general product Bruhat}}. Prove  (a) and (b).
If $w\star w'=ww'$ then the statement 
follows from  (\ref{eq:product Bruhat}). The 
general statement reduces
to the case when $w=w'=s_i$ for $i\in I$. 
Let us use the following 
identity in $SL_2$:
$$\begin{pmatrix} 
c & 0 \\
1 & c^{-1}
\end{pmatrix} \cdot 
\begin{pmatrix} 
c' & 0 \\
1 & c^{'-1}
\end{pmatrix}=
\begin{pmatrix} 
d & 0 \\
1 & d^{-1}
\end{pmatrix} \cdot 
\begin{pmatrix} 
d' & 0 \\
0 & d^{'-1}
\end{pmatrix} \ , 
$$
for any $c,c',d,d'\in \GG_m$ satisfying $d'c=1+dd'$, $dd'=cc'\ne -1$.    
This implies that 
$B^-_{s_i}\cdot B^-_{s_i}\equiv B^-_{s_i}\cdot \alpha_i^\vee(\GG_m)$
and that $B^-_{s_i}\cdot B^-_{s_i}$
contains the set  $B^-_{s_i}\setminus \{\pi_i(-1)\}$, where 
 $\pi_i:\GG_m\widetilde \to B^-_{s_i}$ is the biregular isomorphism
defined in section \ref{subsect:positive unipotent structures}.

Parts (a) and (b) are proved. Parts (c)-(e) easily follow. 

Proposition \ref{pr:general product Bruhat} is proved. 
\hfill$\square$\smallskip

\noindent{\bf Proof of Proposition \ref{pr:torus in product of reduced cells}}.
For  $w\in W$ define an action of $T$ on $G$ by 
\begin{equation}
\label{eq:Adw}
(t,g)\mapsto Ad_w~t(g)=w(t)\cdot g\cdot t^{-1} \ .
\end{equation}

\begin{lemma} 
\label{le:torus in the reduced cell} 
For any $w_1, w_2\in W$ one has 
$Ad_{w_1}~T(U\overline {w_2}U)=U\overline {w_2}U\cdot T_{{w_2}^{-1}w_1}$. 

\end{lemma} 

\begin{proof} It suffices to show that 
$Ad_{w_1}~T(\overline {w_2})=\overline {w_2}\cdot T_{{w_2}^{-1}w_1}$. 
By definition, we have:
$$Ad_{w_1}~t(\overline {w_2})=w_1(t)\cdot \overline {w_2}\cdot t^{-1}=
\overline {w_2}\cdot ({w_2}^{-1}w_1)(t)\cdot t^{-1}
$$
for $t\in T$. Lemma is is proved.
\end{proof}

Furthermore, we have $B^-_w\cdot B^-_{w'}\subset U\overline wU\cdot U\overline {w'}U$. Thus
$Ad_{w\cdot w'} T\left(B^-_w\cdot B^-_{w'}\right)= B^-_w\cdot B^-_{w'}$. On the other hand, 
$Ad_{w\cdot w'} T\left(B^-_{w\star w'}\right)=
B^-_{w\star w'}\cdot T_{(w\star w')^{-1}\cdot w\cdot w'}$.

Therefore, by Proposition \ref{pr:general product Bruhat}(a),   
$B^-_{w\star w'}\cdot T_{(w\star w')^{-1}\cdot w\cdot w'}
\subset B^-_{w{\star}w'}\cdot \widetilde T_{w,w'}$.

Proposition \ref{pr:torus in product of reduced cells} is proved.
\hfill$\square$\smallskip

\noindent{\bf Proof of Proposition \ref{pr:UL-action}}.
Denote temporarily $w:=w_{L',L}$, $b:=\eta^w(x)=\pi^-(x\overline w)$ and 
$u'=Ad~\overline {w_{L',L}}(\pi(u\cdot b))$.  

Furthermore, we have by (\ref{eq:pi for U-action on Uw}): 
$$\pi^w(u\cdot x)= 
\pi\left(u'\cdot \pi^-(\overline w\cdot b^{-1})\right) \ .$$

It is easy to see that  $u'\in U$ because $\pi(u\cdot b)\in U_L$. 
Note that 
$\overline w\cdot b^{-1}\in \overline wU\overline w^{\,-1}U$. 
Thus 
$\pi^-(\overline w\cdot b^{-1})\in
\overline wU_{L'}\overline w^{\,-1}\cap U^-=
\overline w(U_P\cap L')\overline 
w^{\,-1}\subset \overline wU_P\overline 
w^{\,-1}$. The latter set is 
the unipotent radical of $wPw^{\,-1}$. This implies that 
$Ad~u'(\pi^-(\overline w\cdot b^{-1}))\in U^-$ and
$$\pi^w(u\cdot x)=\pi(Ad~u'(\pi^-(\overline w\cdot b^{-1}))\cdot u')=u' \ .$$

Proposition \ref{pr:UL-action} is proved. 
\hfill$\square$\smallskip

\noindent {\bf Proof of Theorem \ref{th:factorization U}}. 
Part (a) follows from Lemma \ref{le:product unipotent bruhat}. Prove (b). 
It suffices to show that $f_{w,w'}$ commutes 
with the action of $U_L$. By definition, the action of $U_L$ 
on the product $({\bf U}^{w_{L',L}},\eta^{w_{L',L}})\times 
({\bf U}^w,\eta^w)$ is given 
by (see \ref{eq:unipotent group action on the product}):
$$u(x,y)=(u(x),\pi\left(u\cdot \eta^{w_{L',L}}(x)\right)(y)) \ .$$

For $x\in U^{w_{L',L}}, y\in U^{w'}$ and $u\in U_L$ we have 
$$u(f_{w,w'}(x,y))=u(x\cdot y')=u(x)\cdot 
u'(y') \ ,$$
where $y'=Ad~ \overline w (y)$ and  $u'=\pi^{w_{L',L}}(u\cdot x)$.
By Proposition \ref{pr:UL-action}, $u'=Ad~\overline {w_{L',L}}(u'')$, 
where $u''=\pi(u\cdot \eta^w(x))$. 

It is easy to see that 
$$u'\cdot y'=
\overline {w_{L',L}}\cdot u''\cdot \overline {w_{L',L}}^{\,-1}\cdot 
\overline w\cdot y\cdot 
\overline w^{\,-1}=Ad~\overline w(u''\cdot y)$$
because $\overline {w_{L',L}}^{\,-1}\cdot 
\overline w$ commutes with $u''\in U_L$.  

Furthermore,
$$u'(y')=u'\cdot y'\cdot \pi^{ww'w^{-1}}(u'y')^{-1}=
Ad~\overline w (u''(y)) $$
because $\pi^{ww'w^{-1}}(u'y')=\pi^{ww'w^{-1}}(Ad~\overline w (u''y))=
Ad~\overline w (\pi^{w'}(u''y)$ 
by Lemma \ref{le:Ad w and pi w} applied with $w''=w_0^L\star w'$.  

Finally,  
$$u(f_{w,w'}(x,y))=u(x)\cdot 
Ad~\overline w (u''(y))=f_{w,w'}\left(u(x),u''(x)\right)(y))=f_{w,w'}(u(x,y))$$ 

Theorem \ref{th:factorization U} is proved. 
\hfill$\square$\smallskip

\noindent{\bf Proof of Theorem \ref{th:chi sigma L invariant}}.
For each $w\in W$ denote $V(w):=U\cap \overline wU^-\overline w^{\,-1}$. 
It is easy to see that $U(w)\cdot V(w)=U$, $U(w)\cap V(w)=\{e\}$, and 
$$U\overline wU=V(w)\overline w U\subset \overline w U^-\cdot U \ .$$

\begin{lemma} 
\label{le:two-sided proper}
For any $w\in W,\chi\in \widehat U$  
we have (see (\ref{eq:chiw})):
\begin{equation}
\label{eq:two-sided proper}
\chi^w(u'gu)=
\overline \chi(\overline w^{\,-1}t^{-1}u't \overline w)+\chi(u)+\chi^w(g)
\end{equation}
for all $g\in t\cdot U\overline w U$, $t\in T$, $u'\in Norm_U(V(w))$ and $u\in U$.  

\end{lemma}

\begin{proof}
Note that $\chi^w(u'gu)=\chi^w(u'g)+\chi(u)$
for all $u$.  So we prove the lemma in the assumption that $u=1$. 

Let $g\in U\overline w U$. Express $g$ as 
$g=u_1 t \overline w u_2$, where $u_1\in V(w)$, $u_2\in U$. 
Then $\chi^w(g)=\overline \chi(\overline w^{\,-1}u_1 t \overline w u_2)=\chi(u_2)$ 
because $\overline w^{\,-1}V(w)\overline w\subset U^-$. Furthermore, 
for each $u'\in Norm_U(V(w))$  we have:
$u'g=u'u_1t\overline wu_2=u'_1u't \overline w  u_2$, where  
$u'_1=Ad~u'(u_1) \in V(w)$. 
Finally,
$$\chi^w(u'g)=\overline \chi(\overline w^{\,-1}u'g)
=\overline \chi(\overline w^{\,-1}u'_1u't\overline w u_2)=
\overline \chi(\overline w^{\,-1}u't\overline w)+\chi(u_2) \ .$$
Lemma is proved. \end{proof}

\begin{proposition} 
\label{pr:chi sigma L invariant}
For any  
$\chi\in \widehat U(w_{L,G})$ and any $t\in Z(L)$ 
the function $f_{\chi,\chi^L}^{w_{L,G}}$ satisfies 
$$f_{\chi,\chi^L}^{w_{L,G}}(u'gu)=\chi(u)+
\chi(u')+f_{\chi,\chi^L}^{w_{L,G}}(g)$$
for $u,u'\in U$, $g\in t\cdot U\overline {w_{G,L}} U$. In particular, 
$f_{\chi,\chi^L}$ is 
invariant under the adjoint action of $U$ on $t\cdot U\overline{w_{G,L}}U$.
\end{proposition}

\begin{proof} Throughout the proof we denote for shortness $w:=w_{L,G}$. 
In this case we have $U(w)=U_L$, $V(w)=U_P$ (where $P=U\cdot B$)  
which implies that $Norm_U(V(w_{L,G}))=U$. 

The following result is obvious.

\begin{lemma} 
\label{le:formula}
Let  $\chi\in \widehat U(w)$. then for any $u\in U$ and $t\in Z(L)$ we have:
$$\overline {\chi^L}\left(\overline {w^{-1}}^{\,-1}\cdot t^{-1}u\cdot 
\overline {w^{-1}}\right)=0,~
\overline \chi(\overline w^{\,-1}u\overline w)=\chi(u)-\chi^L(u) \ .$$
\end{lemma}

Furthermore, let $t\in Z(L)$, $g\in t\cdot U\overline wU$, $u',u\in U$. 
Then 
$$f_{\chi,\chi^L}^{w_{L,G}}(u'gu)=
\chi^w(u'gu)+(\chi^L)^{w^{-1}}(\iota(u)\iota(g)\iota(u'))$$
Applying Lemma \ref{le:two-sided proper} twice,  
we obtain 
$$f_{\chi,\chi^L}^{w_{L,G}}(u'gu)=\overline 
\chi\left(\overline w^{\,-1}t^{-1}u't\overline w\right)+\chi(u)+\chi^w(g)$$
$$+\overline {\chi^L}\left(\overline {w^{-1}}^{\,-1}\cdot t\iota(u)t^{-1}
\cdot \overline {w^{-1}}\right)+
\chi^L(\iota(u'))+(\chi^L)^{w^{-1}}(\iota(g)) \ .$$
By by Lemma \ref{le:formula}, we have
$\overline {\chi^L}(\overline {w^{-1}}^{\,-1}\iota(u) \overline {w^{-1}})=0$, and  
$\overline \chi(\overline w^{\,-1}t^{-1}u't \overline w)=\chi(u')-\chi^L(u')$. 
Taking into account that $\chi^L(\iota(u'))=\chi^L(u')$ we obtain
$$f_{\chi,\chi^L}^{w_{L,G}}(u'gu)=\chi(u')+\chi^w(g)+\chi(u)+\chi^{w^{-1}}(\iota(g)) \ .$$
Proposition \ref{pr:chi sigma L invariant} is proved. 
\end{proof}

Theorem \ref{th:chi sigma L invariant} is proved. 
\hfill$\square$\smallskip

\noindent {\bf Proof of Theorem \ref{th:W invariants}}.
Recall that $W_L=W_J$ denotes the Weyl group of $L=L_J$.

\begin{proposition}  
\label{pr:chi invariant of UL}
For any $\chi\in \widehat U(w_{L,G})$, $t\in T$ we have
$$\chi^{w_{L,G}}(u'gu)=\chi(t^{-1}u't)+\chi(u)+
\chi^{w_{L,G}}(g)$$
for $u'\in U_L, u\in U, g\in t\cdot U\overline {w_{L,G}}U$. 
In particular, $\chi^{w_{L,G}}$ is invariant under the adjoint 
action of $U_L$ on $t\cdot U\overline {w_{L,G}}U$ for $t\in Z(L)$.
\end{proposition}

\begin{proof} Taking into the account that $Norm_U(V(w_{L,G}))=U$, we 
rewrite \ref{eq:two-sided proper}:
$$\chi^{w_{L,G}}(u'gu)=\overline \chi(\overline {w_{L,G}}^{-1}u t
\overline {w_{L,G}})+\chi(u)+\chi^{w_{L,G}}(g)=\chi(t^{-1}u't)+\chi(u)+
\chi^{w_{L,G}}(g) \ .$$

Proposition \ref{pr:chi invariant of UL} is proved. 

\end{proof}

\noindent {\bf Proof of Theorem \ref{th:W invariants}(b)}.
We have  ${\rm supp}(w_{L,G})=I$. By (\ref{eq:simple generator si}), 
for each  $j\in J$ and generic $b_-\in B^-_{w_{L,G}}$ we have 
$s_j(b_-)=Ad~u (b_- )$ for some  $u\in U_i$. This implies that for all
$j\in J$ the restriction 
$\chi^{w_{L,G}}|_{B^-_{w_{L,G}}}$ is  $s_j$-invariant. 
Theorem \ref{th:W invariants}(b) is proved. 
\hfill$\square$\smallskip 

Recall from section \ref{subsect:standard crystals} that the 
$U$-variety ${\bf U}^{w_{G,L}}$ equipped 
with the isomorphism $\eta^{w_{G,L}}:
U^{w_{G,L}}\widetilde \to B^-_{w_{G,L}}$ is an unipotent $G$-crystal. 

\noindent {\bf Proof of Theorem \ref{th:W invariants}(a)}.
The proof is based on the following property of the 
biregular isomorphism  $\eta_w:B^-_w\widetilde \to U^w$ which is defined in  
section \ref{subsect:standard crystals}.

\begin{lemma} For any $\chi\in \widehat U$ and any $w\in W$ we have:
\label{le:chi and eta inverse}
\begin{equation}
\label{eq:chi and eta inverse}
\chi(\eta_w(b))=\chi^{w^{-1}}(\iota(b))
\end{equation}
for $b\in B^-_w$
\end{lemma}

\begin{proof} It is easy to see that for any $\chi\in \widehat U$ and $u\in
U$ we have $\chi(u^{-1})=-\chi(u)$. It is also clear that for any 
$\chi\in \widehat U$, $g\in U^-\cdot T\cdot U$ we have:
$$\overline \chi(t_0gt_0^{-1})=-\overline \chi(g) \ ,$$
where $t_0=\rho^\vee(-1)\in T^{ad}$ as in section \ref{subsect:restriction standard to L}.

Furthermore, by definition of $\overline w$, one has 
$$t_0\overline w t_0^{-1}=\overline {w^{-1}}^{\,-1} \ .$$

Composing $\eta_{w^{-1}}$ with $\overline \chi$, we 
obtain for $b\in B^-_w$, $\chi\in \widehat U$:  
$$ \chi(\eta_w(b))=
\chi(\pi(\overline wb^{-1})^{-1})=
-\chi(\pi(\overline wb^{-1}))
=-\overline\chi(\overline wb^{-1})$$
$$=\overline\chi(t_0\overline wb^{-1}t_0^{-1})=
\chi\left(\overline {w^{-1}}^{\, -1}t_0b^{-1}t_0^{-1}\right)
=\chi^{w^{-1}}(\iota(b)) \ .$$

Lemma is proved. \end{proof}

Now we apply Lemma \ref{le:chi and eta inverse} with 
$w=w_{G,L}=w_{L,G}^{-1}$, $\chi\in \widehat U(w_{G,L})$.  
Using Proposition \ref{pr:chi invariant of UL} and the property of
the mapping $b\mapsto \iota(b)$: 
$$\iota(U\overline {w_{G,L}}U)=U\overline {w_{L,G}}U, 
\iota(Ad~ U_L (g))=Ad~U_L (\iota(g))$$ 
we see that the right hand side of 
(\ref{eq:chi and eta inverse}) is an $Ad~U_L$-invariant function on 
$B^-_{w_{G,L}}$  Therefore, by 
(\ref{eq:simple generator si}), the right hand side of 
(\ref{eq:chi and eta inverse}) is invariant under the 
unipotent action of $W_J$. 

Since, by its definition, $\eta_{w_{G,L}}$ is an isomorphism 
of unipotent $G$-crystals 
$$({\bf B}^-_{w_{G,L}},{\rm id}_{w_{G,L}})\widetilde \to ({\bf U}^{w_{G,L}},
\eta_{w_{G,L}}) \ ,$$
the restriction $\chi|_{U^{w_{G,L}}}$ is also invariant 
under the unipotent action of $W_J$.
   
Theorem \ref{th:W invariants}(a) is proved. 
\hfill$\square$\smallskip

\noindent {\bf Proof of Proposition \ref{pr:factorization B-}}. 
Part (a) follows from (\ref{eq:general product Bruhat}), 
part (b) follows from Corollary \ref{co:bruhat in B-}, part (c) follows 
from Theorem \ref{th:positive unipotent}, and part (d) is clear. 

Proposition \ref{pr:factorization B-} is proved. 
\hfill$\square$\smallskip

\noindent {\bf Proof of Proposition 
\ref{pr:1-dim product unipotent bruhat}}. 
Part (a) is immediate. 
So we prove (b). Due to Lemma \ref{le:the twist U to B-} and 
Proposition \ref{pr:factorization B-}(c), 
it suffices to prove that for any $\ii\in R(w)$ and 
$\ii'\in R(w^{-1})$ the birational isomorphism 
$[\pi_{\ii'}]^{-1} \circ \eta^{w^{-1}}\circ \pi^\ii$ and its inverse are 
positive isomorphisms $T'\widetilde \to T'$, where $T':=(\GG_m)^{l(w)}$. 
This follows from the  results of \cite{bz3}, sections 4 and 5.

Proposition \ref{pr:1-dim product unipotent bruhat} is proved. 
\hfill$\square$

\section {Projections of Bruhat cells to the parabolic subgroups}

\label{sect:projections}

\subsection{General facts on projections}

\label{subsect:general facts}

Recall that we have defined in section 
\ref{subsect:factorization of unipotent $L$-crystals} 
the monoid structure $(W,{\star})$ on $W$.  Clearly, for any standard 
Levi subgroup $L=L_J$, 
the monoid $(W_J,{\star})$ is a sub-monoid of $W$ under the 
operation ${\star}$ defined in section 
\ref{subsect:standard crystals}, and this 
sub-monoid $(W_J,{\star})$ is  
generated by all $s_j$, $j\in J$.

\begin{lemma} The correspondence
$$s_i\mapsto [s_i]=
\begin{cases} 
s_i, & \text{if $i\in J$,} \\
e, & \text{if $i\in I\setminus J$,} \\
\end{cases}$$
extends to the  homomorphism  of monoids $[\cdot]:
(W,{\star})\to (W_J,{\star})$. 
\end{lemma}

Similarly to the projection ${\bf p}_L^-:B^-\to B^-_L$ defined in section 
\ref{subsect:L-crystals}
let  ${\bf p}^+={\bf p}_L^+$ be the natural projection $U\to U_L$

There is a natural equivalence relations between constructible  
subsets of $G$. We say that two such subsets $S,S'\subset G$ 
are {\it equivalent} if
their intersection $S\cap S'$ is dense in each $S$ and $S'$. 
In this case we denote $S\equiv S'$.

\begin{lemma} 
\label{le:dense subset}
For each $w\in W$ one has ${\bf p}^+(U^w)\equiv U^{[w]}$.

\end{lemma} 

\begin{proof} We proceed by the induction in $l(w)$. 
If $l(w)=1$, that is, 
$w=s_i$ for some $i\in I$ then
\begin{equation}
\label{eq:simple projection}
{\bf p}^+(U^{s_i})=
\begin{cases} U^{s_i},  & \text{if $i\in J$,} \\
              \{e\},  & \text{if $i\in I\setminus J$.}
\end{cases} 
\end{equation}

Let us now assume that $l(w)>1$. Using Proposition
\ref{pr:1-dim product unipotent bruhat} and the fact that ${\bf p}^+$ 
is a group homomorphism, we 
obtain for any $\ii=(i_1,\ldots,i_l)\in R(w)$:
$${\bf p}^+(U^w)\equiv {\bf p}^+(U^{s_{i_l}}U^{s_{i_{l-1}}}\cdot 
\cdots \cdots U^{s_{i_1}})$$
$$={\bf p}^+(U^{s_{i_l}})\cdot {\bf p}^+(U^{s_{i_{l-1}}})\cdot 
\cdots \cdots {\bf p}^+(U^{s_{i_1}})\equiv 
U^{[s_{i_1}]{\star}\cdots {\star}[s_{i_l}]}\equiv U^{[w]}$$

Lemma \ref{le:dense subset} is proved. \end{proof}

The computation of ${\bf p}_L^-(B^-_w)$ is a little more complicated.

\begin{proposition}
\label{pr:reduced whittaker} For any $w\in W$ we have:

\noindent (a) The intersection of ${\bf p}_L^-(B^-_w)$ with 
$B^-_{[w]}$ is a dense subset of $B^-_{[w]}$;

\noindent (b) There is an algebraic sub-torus 
$\widetilde T_w\subset T$ such that ${\bf p}_L^-(B^-_w)\equiv B^-_w\cdot \widetilde T_w$.

\noindent (c) The torus $\widetilde T_w$ satisfies: $\widetilde T_w=\{e\}$ if and only if 
$w\in W_J$, $\widetilde T_{s_i}=\alpha_i^\vee(\GG_m)$ for $i\in I\setminus J$, and 
$\widetilde T_{ww'}= [w']^{-1}(\widetilde T_w)\cdot \widetilde T_{w'}\cdot 
\widetilde T_{[w],[w']}$
for any $w,w'$ such that $w\star w'=ww'$. 

\end{proposition}

\begin{proof} It follows from Proposition 
\ref{pr:general product Bruhat} that for any $w,w'\in W$ 
satisfying $w\star w'=ww'$ we have:
$${\bf p}_L^-(B^-_{ww'})\equiv {\bf p}_L^-(B^-_w\cdot B^-_{w'})= 
{\bf p}_L^-(B^-_w)\cdot {\bf p}_L^-(B^-_{w'})\equiv B^-_{[w]}\cdot \widetilde T_w\cdot 
B^-_{[w']}\cdot \widetilde T_{w'}$$
$$\equiv B^-_{[w]}\cdot B^-_{[w']}\cdot [w']^{-1}(\widetilde T_w)\cdot 
\widetilde T_{w'}\cdot \widetilde T_{[w],[w']} \ .$$

This proves (b) and (c). Let us prove (a) now. Denote 
$\tilde B^-_{[w]}:={\bf p}_L^-(B^-_w)\cap B^-_{[w]}$. We will proceed by 
induction in $l(w)$. If $w=e$, the statement is obvious. Let now $w\ne e$. Let us express 
$w=w'w''$, where $w'\ne e$, $w''\ne e$ and $l(w)=l(w')+l(w'')$. Thus, the 
inductive hypothesis implies that $\tilde B^-_{[w']}$ is dense in $B^-_{[w']}$, 
and $\tilde B^-_{[w'']}:={\bf p}_L^-(B^-_{w''})\cap B^-_{[w'']}$ is 
dense in $B^-_{[w'']}$.  

Then (\ref{eq:product Bruhat}) implies that 
$${\bf p}_L^-(B^-_{w'w''})\supset {\bf p}_L^-(B^-_{w'}\cdot B^-_{w''})\supset 
{\bf p}_L^-(B^-_{w'})\cdot 
{\bf p}_L^-(B^-_{w''})\supset \tilde B^-_{[w']}\cdot \tilde B^-_{[w'']}$$
But $\tilde B^-_{[w']}\cdot \tilde B^-_{[w'']}$ is a 
dense subset of $B^-_{[w']}\cdot B^-_{[w'']}$. Finally, Proposition 
\ref{pr:general product Bruhat}(a)  implies that the latter set 
contains a dense subset of the variety $B^-_{[w']{\star}[w'']}=B^-_{[w'w'']}$. 

This proves part (a). Proposition is proved. \end{proof}

It is well-known that the set $G_0=U^-_P\cdot P$ is open in $G$. 
Denote by  ${\bf p}={\bf p}_L$ the natural projection $G_0\to P$. By definition, 
${\bf p}$ commutes with the action of  $B_L\times B$, and 
the restriction of ${\bf p}$ to $B^-\subset G_0$ is the 
natural projection ${\bf p}_L^-:B^-\to B^-_L$. 

Our next task is to describe the image of each reduced Bruhat cell 
$U\overline wU$ under the projection ${\bf p}$.

Recall that in section 
\ref{subsect:factorization of unipotent $L$-crystals} we defined 
for each $w\in W$ the sub-torus $T_w\subset T$.

\begin{proposition} 
\label{pr:whittaker} 
For any $w\in W$ we have:

\noindent (a) $U\overline {[w]}U\cdot T_{[w]^{-1}w}\subset 
{\bf p}(U\overline w U\cap G_0)$. 
 
\noindent (b) The restriction ${\bf p}_w={\bf p}_{L;w}$  of ${\bf p}$ 
to $U\overline wU$ is a 
dominant rational morphism 
$${\bf p}_w:U\overline wU \to U\overline {[w]}U\cdot \widetilde T_w \ .$$

\noindent (c) $T_{[w]^{-1}w}\subset \widetilde T_w$. 
\end{proposition}

\begin{proof} Prove (b). 
The multiplication in $G$ induces the open inclusions 
$B^-_w\times U\hookrightarrow U\overline wU\cap G_0$, 
$B^-_w\times U_L\hookrightarrow U\overline wU_L\cap G_0$. Furthermore,
\begin{equation}
\label{eq:double projection to single}
{\bf p}(B^-_w\cdot U)={\bf p}_L^-(B^-_w)\cdot U,~ 
{\bf p}(B^-_w\cdot U_L)={\bf p}_L^-(B^-_w)\cdot U_L, 
\end{equation}

Therefore, it follows from Lemma \ref{le:dense subset} that
$${\bf p}(U\overline wU)\equiv  
{\bf p}_L^-(B^-_w)\cdot U\equiv  
B^-_{[w]}\cdot \widetilde T_w\cdot U=B^-_{[w]}\cdot U 
\cdot \widetilde T_w\equiv U\overline 
{[w]}U\cdot \widetilde T_w \ .$$ 

Part (b) is proved. 

Let us prove (a) now. 
Recall that ${\bf p}$ is  $U_L\times U$-equivariant, and
$U\overline {[w]}U=U_L\overline {[w]}U$. Note that 
${\bf p}_L^-(B^-_w)$ intersects $B^-_{[w]}$ non-trivially. Thus 
$$U\overline {[w]}U \subset {\bf p}(U\overline wU\cap G_0)$$

Recall that the action $Ad_w:T\times G\to G$ is 
defined in the proof of Proposition 
\ref{pr:torus in product of reduced cells}. 
Clearly, both $P$ and the reduced  Bruhat cell 
$U\overline wU$ are invariant under this action, and the morphism  
${\bf p}:G_0\to P$ commutes with this action. Thus applying  Lemma 
\ref{le:torus in the reduced cell} to the above identity, we obtain 
$$Ad_w~T(U\overline {[w]}U)=U\overline {[w]}U\cdot T_{{[w]}^{-1}w}
\subset {\bf p}(U\overline wU\cap G_0) \ .$$

Part (a) is proved. Part (c) follows. 

Proposition \ref{pr:whittaker} is proved. \end{proof}

Denote  $U^-_{P;w}:=U_P^-\cap \left(U\overline w U\cdot \overline{[w]}^{\,-1}\right)$. 
The set $U^-_{P;w}$ is not empty because 
$\overline {[w]}\in {\bf p}(U\overline w U\cap G_0)$ 
due to Proposition \ref{pr:whittaker}(a). 
Note also that $U^-_{P;w}$ is closed in $U^-_P$. 
For any $u_-\in U^-_{P;w}$ let 
${\bf q}_{u_-}: U\overline {[w]}U\cdot T_{[w]^{-1}w}\to U\overline w U$ 
be the morphism defined by 
$${\bf q}_{u_-}\left (u'\cdot \overline {[w]}\cdot u\cdot t\right):=
u'\cdot u_-\cdot \overline {[w]}\cdot u\cdot t$$
for $u\in U, u'\in V([w])$, 
and $t\in T_{[w]^{-1}w}$.

For any $w\in W$ let $T^w:=\{t\in T:w(t)=t\}$. Clearly, $\dim~T^w\cdot
T_w=\dim~T$.

\begin{proposition} 
\label{pr:qq}
For any $w\in W$ we have:

\noindent (a) for each  $u_-\in U^-_{P;w}$ the morphism 
${\bf q}_{u_-}$ is a section of ${\bf p}_w$, that is, 
$${\bf p}_w\circ {\bf q}_{u_-}={\rm id}_{U\overline 
{[w]}U\cdot T_{[w]^{-1}w}} \ . $$

\noindent (b) The variety $U^-_{P;w}$ is invariant under the 
adjoint action of $T^{[w]^{-1}\cdot w}\cdot U_L([w])$. 

\noindent (c) Let ${\bf q}$ be a morphism $U^-_{P;w}\times 
U\overline {[w]}U\cdot T_{[w]^{-1}w}\to U\overline w U$ defined by 
${\bf q}(u_-,g)={\bf q}_{u_-}(g)$. Then ${\bf q}$ is an inclusion.
\end{proposition}

\begin{proof} Prove (a). Indeed, $V([w])\subset U_L$. Thus ${\bf q}_{u_-}$
is a section of ${\bf p}_w$ because ${\bf p}_w$ is 
$U_L\times B$-equivariant. 

Prove (b). It is easy to see that $U_L([w])=Norm_{U_L}(\overline 
{[w]})$. Thus, for $u_-\in U^-_{P;w}$ and any $u\in Norm_{U_L}(\overline
{[w]})$, 
$t\in T^{[w]^{-1}\cdot w}$ we have 
$Ad_w~t(u\cdot u_-\cdot \overline {[w]}\cdot u^{-1})\in U\overline wU$. On
the other hand, 
$Ad_w~t(u\cdot u_-\cdot \overline {[w]}\cdot u^{-1})=
(tu)\cdot u_-\cdot (tu)^{-1}\cdot \overline {[w]}$.  This implies that 
$(tu)\cdot u_-\cdot (tu)^{-1}\in U^-_{P;w}$.

Prove (c). We proceed by the contradiction.  Assume that ${\bf q}$ is not
injective. 
That is, we assume that there are elements $u_-,u'_-\in U^-_{P;w}$ such that 
$u'_-\ne u_-$ and 
$${\bf q}_{u'_-}(\overline {[w]})={\bf q}_{u_-}(u'\overline {[w]}ut)$$
for some $u'\in V([w]), u\in U$ and $t\in T$. In other words, 
$$u'_-\cdot \overline {[w]}=u'\cdot u_-\cdot 
\overline {[w]}\cdot u\cdot t \ .$$
Denote $\tilde u_-:= Ad~u'(u_-)$. 
Clearly, $\tilde u_-\in U_P^-$ because $u'\in U_L$.  Furthermore, let us
factorize
$Ad~\overline {[w]}(u)=u_+\cdot u''_-$, where $u_+\in U(\overline {[w]})$, and 
$u''_-\in U_L^-$. Thus, the above identity can be rewritten as 
$$(\tilde u_-)^{-1}\cdot u'_-= u'\cdot u_+\cdot u''_-\cdot [w](t) 
\ .$$

This implies that $t=e$, $u''_-=e$, $u'\cdot u_+=e$, 
and $u'_-=\tilde u_-$. In particular, $u'=(u_+)^{-1}\in V([w])\cap
U([w])=U_L([w])=\{e\}$. 

In its turn, this implies that $u'_-= u_-$ which contradicts to the
original assumption.

\end{proof}

\noindent {\bf Definition}. We call an element $w\in W$ 
{\it $L$-special} (or simply {\it special}) if
$l(w)=l([w])+\dim~T_{[w]^{-1}w}$.

\begin{theorem} 
\label{th:special isomorphism}
For any special element $w\in W$ we have:

        \noindent (a) The set $U_{P;w}^-$ consists of a single element
        $u_-=u^-_{P;w}$. This element $u_-$ is centralized by the group 
        $T^{[w]^{-1}\cdot w}\cdot U_L([w])$. The inclusion 
        ${\bf q}_{u_-}$
        is dense.

\noindent (b) ${\bf p}_w=[{\bf q}_{u_-}]^{-1}$. In particular, 
${\bf p}_w:U\overline w U \to U\overline{[w]} U\cdot T_{[w]^{-1}w}$ 
is a birational isomorphism.

\noindent (c) $\widetilde T_w= T_{[w]^{-1}w}$.

\noindent (d) The composition ${\bf p}^-_L\circ \eta^w$ is a birational
isomorphism $U^w\widetilde \to B^-_{[w]}\cdot T_{[w]^{-1}w}$. 

\end{theorem}

\begin{proof} Follows from Proposition \ref{pr:qq}. 

\end{proof}

\subsection {Example of the projection for $G=GL_{m+n}$}

Let us fix two positive integers $m\le n$.  
Let $G=GL_{m+n}$ and $L=L_{m,n}=GL_m\times GL_n\subset GL_{m+n}$. 
In this case 
$P=P_{m,n}=L\cdot B$ is the corresponding maximal parabolic subgroup 
which consists of 
those matrices $g\in GL_{m+n}$ which have zeros in the 
$m\times n$-rectangle in the lower left corner. 
  
We use the standard labeling of the Dynkin diagram of type $A_{m+n-1}$: 
$I=\{1,2,\ldots,m+n-1\}$. Let $J=I\setminus \{m\}$.  Then the Weyl group $W$ is naturally 
identified with the symmetric group $S_{m+n}$ via $s_i\mapsto (i,i+1)$.

Note that the longest element 
$w_0$ of $W$ acts on simple roots by $w_0(\alpha_i)=-\alpha_{i^*}$ for 
$i=1,\ldots,m+n-1$, where  $i^*=m+n-i$. 

The corresponding element  $w_{L,G}=w_{m,n}\in W$ admits a 
factorization (\ref{eq:wmn}). 

For $i\le j$ denote by  $s_{ij}$ the reflection about the root 
$\alpha_i+\cdots +\alpha_j$. Clearly, $s_{ij}$ is the 
transposition $(i,j+1)\in W=S_{m+n}$.

Denote $\sigma_{m,n}:=[w_{m,n}]^{-1}w_{m,n}$.
It is easy to see that 
\begin{equation}
\label{eq:commuting sij}
\sigma_{m,n}=\prod_{i=1}^m s_{i,i^*} \ ,
\end{equation}
that is, 
$\sigma_{m,n}$  is a product of exactly $m$ 
commuting reflections in $W$. This implies that 
$\dim~T_{\sigma_{m,n}}=m=l(w_{m,n})-l([w_{m,n}])$. 
Thus $w_{m,n}$ is special. 

Let ${\bf q}_{m,n}$ be the dense inclusion  
$U_L\overline{[w_{m,n}]}U\cdot T_{\sigma_{m,n}}
\hookrightarrow U\overline {w_{m,n}}U$ 
prescribed by Theorem \ref{th:special isomorphism}.

Let $\omega_1,\ldots,\omega_{m+n}$ be the fundamental 
weights for $GL_{m+n}$. For each $w\in S_{m+n}$ 
the weight $w(\omega_k)$ is naturally identified with a $k$-element subset 
of $\{1,\ldots,m+n\}$. 

For any $w_1,w_2\in S_{m+n}$ let 
$\Delta_{w_1(\omega_i),w_2(\omega_i)}:G\to \AA^1$ 
be the minor located in the intersection of rows 
indexed by $w_1(\omega_i)$ and columns indexed by 
$w_2(\omega_i)$.

\begin{proposition}
\label{def:M-}
The variety $U_L\overline{[w_{m,n}]}U_L\cdot T_{\sigma_{m,n}}$ consists of all 
$g\in B_L\overline{[w_{m,n}]}B_L$ satisfying for $i=1,\ldots,m+n-1$:
\begin{equation}
\label{eq:equality of Deltas1}
\Delta_{[w_{m,n}](\omega_i),\omega_i}(g)=
\begin{cases} 
\Delta_{\omega_m,\omega_m}(g), & \text{if $i\in [m;n]$,} \\
\Delta_{[w_{m,n}](\omega_i^*),\omega_i^*}(g), & \text{if $i\notin [m;n]$.} 
\end{cases} 
\end{equation}

\end{proposition}

\begin{proof} Let us describe the sub-torus $T_{\sigma_{m,n}}$ of $T$. We have 
$$T_{\sigma_{m,n}}=\prod_{i=1}^m T_{s_{i,i^*}} \ .$$
That is, $T_{\sigma_{m,n}}$ can be described in $T$ by the equations 
($i=1,\ldots,m+n-1$): 
$$\Delta_{\omega_i,\omega_i}(t)=
\begin{cases} 
\Delta_{\omega_m,\omega_m}(t), & \text{if $i\in [m;n]$,} \\
\Delta_{\omega_i^*,\omega_i^*}(t), & \text{if $i\in [m;n]$.} 
\end{cases} $$

On the other hand, for any $w\in W_J, t\in T$ one has 
$$U_L\overline w U_Lt=\{g\in B_LwB_L:
\Delta_{w(\omega_i),\omega_i}(g)=\Delta_{\omega_i,\omega_i}(t),i=1,\ldots,
m+n-1\} \ .$$

Proposition \ref{def:M-} is proved. 
\end{proof}

\begin{proposition} For any $\chi\in \widehat U(w_{m,n})$,  
$u'\in U_L, u\in U$, $t\in T$ we have:
\begin{equation}
\label{eq:two sided proper function pchiL}
\chi^{w_{m,n}}({\bf q}_{m,n}(u' \overline {[w_{m,n}]} \cdot 
\sigma_{m,n}(t)\cdot t^{-1}\cdot u)
=\chi(u')+\chi(u)-\delta_{m,n}\alpha_m(t)\cdot\chi(x_m(1)) \ .
\end{equation}
\end{proposition} 

\begin{proof} 
For any $1\le i\le m\le j$, $w\in W_J$ let $\tilde s_{ij}\in Norm_G(T)$ be the 
representative of $s_{ij}$ defined by 
\begin{equation}
\label{eq:tilde s}
\widetilde s_{ij}:=
\overline w^{\, -1}\overline s_m \overline w
\end{equation}
for any $w\in W_J$ such that $w^{-1} s_mw=s_{ij}$. 

We need the following obvious refinement of (\ref{eq:commuting sij}). 

\begin{lemma}
\label{le:commuting tilde sij}
$\overline{[w_{m,n}]}^{\,-1}\overline w_{m,n}=\prod_{i=1}^m 
\widetilde {s_{i,i^*}}$.
\end{lemma}

For any  $1\le i\le m$ choose an element $w_i\in W$ such that 
$s_{i,i^*}=w_i^{-1}\overline s_mw_i$. By 
(\ref{eq:split si}), we have: 
$$\widetilde {s_{i,i^*}}^{\,-1}=u_i^+u_i^-u_i^+, $$ 
where $u_i^+=\overline w_i^{\, -1}x_m(-1) \overline w_i,u_i^-=
\overline w_i^{\, -1}y_m(1) \overline w_i$. 

Clearly,  $u_i^+\in U_P\cap 
Ad~\overline {w_{m,n}}^{\,-1}(U^-_P)$, and  
$u_i^-\in U^-_P\cap Ad~\overline {w_{m,n}}^{\,-1}(U_P)$. It is also clear 
that  $u_i^\varepsilon$ commutes with each 
$u_j^{\varepsilon'}$ whenever $i\ne j$, where $\varepsilon,\varepsilon'\in \{+,-\}$. 

Denote 
$u^\varepsilon:=\prod_{i=1}^m u_i^\varepsilon$ for 
$\varepsilon\in \{+,-\}$.  

\begin{lemma} 
\label{le:Up}
The only element of $U_{P,w_{m,n}}$ is equal to 
$Ad~\overline{w_{m,n}}((u^+)^{-1})$, and 
$$u_-\cdot\overline {[w_{m,n}]}=\overline{w_{m,n}} \cdot u^-\cdot u^+ \ .$$ 
\end{lemma}
 
\begin{proof} Let $u_-:=Ad~\overline{w_{m,n}}((u^+)^{-1})$. Clearly, $u_-\in U_P^-$. 

Next, let us rewrite 
the equation from Lemma \ref{le:commuting tilde sij}:
$\overline{w_{m,n}}^{\,-1}\overline {[w_{m,n}]}=u^+\cdot u^-\cdot u^+$
or, equivalently, 
$\overline{[w_{m,n}]}=\overline {w_{m,n}}\cdot
u^+\cdot u^-\cdot u^+
=(u_-)^{-1}\cdot u_1^+ \overline {w_{m,n}}\cdot u^+$, where 
$u_1^+=Ad~\overline{w_{m,n}} (u^-)$. Clearly, $u_1^+\in U_P$.  

Thus,  $u_-\in U_P^-\cap U\overline{w_{m,n}}U\overline{[w_{m,n}]}^{\,-1}=U_{P,w_{m,n}}$.  

Lemma \ref{le:Up} is proved. \end{proof}

Since ${\bf q}_{m,n}$ is $U_L\times U$-equivariant it suffices to prove  
(\ref{eq:two sided proper function pchiL}) in the case when $u=u'=e$. 
Using Lemma \ref{le:Up} we obtain: 
$$\chi^{w_{m,n}}({\bf q}_{m,n}(Ad_{w_{m,n}}~t(\overline {[w_{m,n}]}))=
\chi^{w_{m,n}}(Ad_{w_{m,n}}~t({\bf q}_{m,n}(\overline {[w_{m,n}]}))$$
$$\chi^{w_{m,n}}(Ad_{w_{m,n}}~t(\overline {w_{m,n}}u^-u^+))
=\overline \chi(u^-u^+\cdot t^{-1})=\chi(tu^+t^{-1})\ .$$

Let us compute $\chi(tu^+t^{-1})$ using the fact that 
$u_i^+\in [U,U]$ for $i=1,\ldots,m-1$: 
$$\chi(tu^+t^{-1})=\sum_{i=1}^m \chi(tu_i^+t^{-1})=\chi(tu_m^+t^{-1})=
-\delta_{m,n}\chi(x_m(t))$$
The last equality holds because because $tu_m^+t^{-1}\in [U,U]$ for $m\ne n$, and  
$tu_m^+t^{-1}=tx_m(-1)t^{-1}=x_m(-\alpha_m(t))$ if $m=n$.

Proposition is proved. \end{proof}

Let $L'=L'_{n-m}:=(\GG_m)^{2m}\times GL_{n-m}\subset L_{m,n}$. By definition, 
$[w_{m,n}]=w_{{L',L_{m,n}}}$. It is easy to see 
that $[w_{m,n}](T_{\sigma_{m,n}})\subset Z(L')$. 

Recall from Proposition \ref{pr:chi sigma L invariant} that for any 
$\chi\in \widehat U$ the function 
$f_{\chi,\chi^{L'}}^{[w_{m,n}]}$ on $Z(L')\cdot U_L\overline {[w_{m,n}]}U_L$ defined by 
(\ref{eq:fchi}) is proper under the action of $U_L\times U_L$. 

\begin{corollary} For any $\chi\in \widehat U(w_{m,n})\cap \widehat U([w_{m,n}])$, and 
$g\in U_L\overline{[w_{m,n}]} U_L\cdot t$, $t\in T_{\sigma_{m,n}}$  one has:
$$\chi^{w_{m,n}}({\bf q}_{m,n}(g))
=f_{\chi,\chi^{L'}}^{[w_{m,n}]}(g)-\delta_{m,n}a_m\cdot t_m^{-1} \ .$$
where $a_m$ is the coefficient of $\chi_m$ in the expansion 
$\chi=\sum_{i=1}^{m+n-1} a_i \chi_i$.
\end{corollary}

\noindent {\bf Remark}.  It is well-known that $U_{P_{m,n}}$ is 
abelian in the above example. One can expect that for any   
$G$ and any standard parabolic subgroup $P=L\cdot U_P$ the element $w_{L,G}$ 
is special if and only if  $U_P$ is abelian. 

\smallskip  

On the other hand, there is an example of a standard parabolic subgroup in the simple 
group $G$ of the type $D_4$ with non-abelian $U_P$ and non-special  
$w_{L,G}$.

We use the following labeling of the Dynkin diagram of type $D_4$: 
$I=\{0,1,2,3\}$, where $0$ stands for the ``center" of the Dynkin diagram. 
Let $J=\{1,2,3\}$, and $L:=L_J$.  
Then $w^L_0=s_1s_2s_3$, and $w_{L,G}=s_0s_1s_2s_3s_0s_1s_2s_3s_0$.  
It is easy to see that $[w_{L,G}]=w^L_0$. Thus,
$l(w_{L,G})-l([w_{L,G}])=9-3=6>dim~T_{[w_{L,G}]^{-1}w_{L,G}}$, 
which implies 
that $w_{L,G}$ is not special. 

Moreover,  ${\bf p}_w$ is not a birational isomorphism because 
 $\dim ~U\overline {w_{L,G}}U=21$, $\dim~ U\overline {[w_{L,G}]}U= 15$, 
and $\dim~ \widetilde T_{w_{L,G}}\le 4$.

\section{Appendix: Combinatorial pre-crystals and $W$-crystals}

\label{append:pre-crystals}

In this section we recall  
Kashiwara's framework (see \cite{k93}). 

\noindent {\bf Definition}. A {\it partial bijection} of sets 
$f:A\to B$ is a bijection $A'\to B'$ of subsets $A'\subset A$, 
$B'\subset B$. We denote the subset $A'$ by 
$dom(f)$ and  the subset $B'$ by $ran(f)$.

\smallskip

\noindent {\bf Remark}. A partial bijection  $f:A\to B$ is an 
embedding $A\hookrightarrow B$ if and only if  $dom(f)=A$.

\smallskip

The inverse $f^{-1}$ of a partial bijection $f:A\to B$ is the inverse 
bijection $ran(f)\to dom(f)$. 
Composition  $g\circ f$ of partial bijections $f:A\to B, g:B\to C$ 
is naturally a partial bijection with 
$dom(g\circ f)=dom(f)\cap f^{-1}(ran(f)\cap dom(g))$ and 
$ran(g\circ f)=g(ran(f)\cap dom(g))$. In particular, for 
any partial bijection 
$f:B\to B$ and $n\in \ZZ$ the  $n$-th power $f^n$ is a 
partial bijection $B\to B$.

Note that for any partial bijection $f:A\to B$ 
the composition $f^{-1}\circ f$ is the 
{\it partial identity bijection} $id_{dom(f)}:A\to A$.  
   
\bigskip

\noindent {\bf Definition}. Let $B$ be a set 
$\tilde \gamma:B\to \Lambda^\vee$ be a map. We call a family  
of partial bijections {\it combinatorial pre-crystal} (or simply 
{\it pre-crystal}) on $(B,\tilde \gamma)$ if
$$\tilde \gamma(\tilde e_i(b))=\tilde \gamma(b)+\alpha_i^\vee$$ 
for $b\in dom(\tilde e_i)$, $i\in I$.  

\smallskip

With any pre-crystal ${\mathcal B}$  we define partial bijections 
$\tilde s_i:B\to B$, $i\in I$, by  
\begin{equation}
\label{eq:Kashiwara W}
\tilde s_i(b)=\tilde e_i^{\,-\left<\tilde \gamma(b),\alpha_i\right>}(b) \ .
\end{equation}
In fact  $\tilde s_i$ are partial involutions: 
$\tilde s_i^2=id_{dom(\tilde s_i)}$. 

\medskip

\noindent {\bf Examples}.

\noindent 1. Fix $i\in I$. Denote by  $\tilde \gamma_i:\ZZ\to \Lambda^\vee$ 
the map $n\mapsto n\alpha_i^\vee$. We define a pre-crystal  ${\mathcal
B}_i$ on $(\ZZ ,\tilde \gamma_i)$ in the following way. We define a
partial bijection $\tilde e_i:\ZZ\to \ZZ$ by  
$\tilde e_i(n)=n+1$ and for any $j\ne i$ we define
$dom(\tilde e_j)=ran(\tilde e_j)=\emptyset$. 
This defines a pre-crystal on $(\ZZ,\tilde \gamma_i)$ which we denote by
${\mathcal B}_i$ (see \cite{k93}, Example 1.2.6). 

\smallskip

\noindent 2. The lattice $\Lambda^\vee$ is a pre-crystal with 
$\tilde \gamma={\rm id}|_{\Lambda^\vee}$ and 
$\tilde e_i(\lambda^\vee)=\lambda^\vee+\alpha_i^\vee$
for $\lambda^\vee\in \Lambda^\vee$, $i\in I$. 

\smallskip

\noindent 3. For any pre-crystal ${\mathcal B}$ on $(B,\tilde \gamma)$ 
we denote by
$\tilde \gamma'$ the map $\Lambda^\vee\times B \to 
\Lambda^\vee$ given by $\tilde \gamma'(\lambda^\vee,b)=
\lambda^\vee+\tilde \gamma(b)$. 
One can define a pre-crystal on  
$(\Lambda^\vee \times B,\tilde \gamma_{\Lambda^\vee})$  by 
$\tilde e_i(\lambda^\vee, b)=(\lambda^\vee,\tilde e_i(b))$. 

\smallskip

\noindent 4. For any pre-crystal ${\mathcal B}$ on $(B,\tilde \gamma)$ put 
$\gamma^*:=-\gamma$. Define also the partial bijections
$\tilde e_i^{\,*}:=\tilde e_i^{\,-1}$. The collection $\tilde e_i^{\,*}$, 
$i\in I$, defines 
the structure of a pre-crystal on $(B,\gamma^*)$. 
This pre-crystal is called {\it dual} to ${\mathcal B}$ and is denoted by 
${\mathcal B}^*$.

\medskip

\noindent {\bf Remarks}. 

\noindent 1. In the \cite{k93} the partial bijection 
$\tilde e_i^{\,-1}$ was 
denoted by $\tilde f_i$.

\noindent 2. Our pre-crystals correspond to Kashiwara's 
crystals for ${}^L\gg$.

\smallskip

\noindent {\bf Definition}.
A pre-crystal ${\mathcal B}$ is called {\it free} if each $\tilde e_i$ 
is a bijection $B\to B$. 

\smallskip

\noindent {\bf Definition}. A {\it morphism} of pre-crystals
$f:{\mathcal B}'\to {\mathcal B}$ is  a partial
bijection $f:B'\to B$  such that $dom(f)=B',\tilde \gamma=\tilde
\gamma'\circ f$, and the partial bijections $\tilde
e'_i:B'\to B'$, $i\in I$, are obtained from the partial bijections
$e_i:B\to B$, $i\in I$, by $\tilde e'_i=f^{-1} \circ \tilde e_i\circ f$. 

\smallskip

\noindent {\bf Definition}. We say that a pre-crystal ${\mathcal B}$ is 
a {\it combinatorial $W$-crystal} (or simply $W$-{\it crystal})  
if for any $i\in I$ we have: $dom(\tilde s_i)=B$,   
and the involutions $\tilde s_i$ satisfy the braid relations.

It is clear that a structure of a $W$-crystal on  $(B,\tilde \gamma)$
defines an action of $W$ on $B$ in such a way that for any $i\in I
~s_i\in W$ acts by $\tilde s_i:B\to B$.

\smallskip

\noindent {\bf Remarks}.

\noindent 1. The condition that $dom(\tilde s_i)=B$ 
is equivalent to:
\begin{equation}
\label{eq:domain of si}
\inf\{n:b\in dom(\tilde e_i^{\,n})\}\le 
-\left<\tilde \gamma(b),\alpha_i\right>  \le \sup\{n:b\in dom(\tilde e_i^{\,n})\}
\end{equation} 
for all $b\in B, i\in I$.

\noindent 2. The pre-crystal on $(\Lambda^\vee,{\rm id}_{\Lambda^\vee})$ 
defined above  is a free $W$-crystal. 

\noindent 3. For any $W$-crystal the structure map 
$\tilde \gamma:B\to \Lambda^\vee$ 
is $W$-equivariant. 

\medskip

 We denote by  $Pre-Cryst$ the category such that the
objects are pre-crystals and arrows are morphisms of pre-crystals. 
We denote by $W-Cryst$ the full sub-category of $Pre-Cryst$ whose 
objects are $W$-crystals. Note that each morphism 
$f:{\mathcal B}\to {\mathcal B}'$ 
in $W-Cryst$ is an injective $W$-equivariant map $B\to B'$.

\smallskip

\noindent  {\bf Example}. Let $V$ be a finite-dimensional 
$U_q({}^L\gg)$-module, and 
let $B=B(V)$ be a  
{\it crystal basis} for $V$ (see \cite{k90}). Denote by 
$\tilde \gamma:B\to \Lambda^\vee$ 
the weight grading. Then the partial bijections 
$\tilde e_i,\tilde f_i:B\to B$
define a structure of a pre-crystal on $(B,\gamma)$. 
It is a deep result of \cite{k94} that  this is a $W$-crystal.

We denote this $W$-crystal by ${\mathcal B}(V)$.  

\medskip
  
\noindent {\bf Remark}.  
For ${\mathcal B}={\mathcal B}(V)$ all $e_i,e_i^{-1}$ are nilpotent, that is, 
$dom(\tilde e_i^{\,n})=dom(\tilde e_i^{-\,n})=\emptyset$ for some $n>0$,  and 
$-\left<\tilde \gamma(b),\alpha_i\right>=\inf\{n:b\in dom(\tilde e_i^{\,n})\}+ 
\sup\{n:b\in dom(\tilde e_i^{\,n})\}$ 
for $b\in {\mathcal B}(V)$, $i\in I$. Clearly, this identity implies 
the inequalities (\ref{eq:domain of si}), that is, $dom(\tilde s_i)={\mathcal B}(V)$ 
for each $i\in I$.

\smallskip

There are examples of free $W$-crystals, i.e., in which each $\tilde e_i$ is torsion-free. 
Among them are the free $W$-crystals ${\mathcal B}_\ii$ 
from  \cite{k93}, section 2.2 built on $(\ZZ^{l(w_0)},\tilde \gamma_\ii)$, 
$\ii\in R(w_0)$. Kashiwara has constructed each ${\mathcal B}_\ii$ as 
the {\it tensor product} of ${\mathcal B}_{i_1},\ldots,{\mathcal B}_{i_l}$, 
so that ${\mathcal B}_\ii$ is equipped with the family of functions 
$\tilde \varphi^{\ii}_i:\ZZ^{l(w_0)}\to \ZZ$, $i\in I$, 
satisfying
$\tilde \varphi^{\ii}_i(\tilde e_i(b))
=\tilde \varphi^{\ii}_i(b)+1$ 
for $b\in \ZZ^{l(w_0)}$, $i\in I$.

Let ${\mathcal R}(w_0)$ be the category 
whose objects are reduced sequences $\ii\in R(w_0)$ and 
for each pair of objects
$\ii,\ii'\in R(w_0)$ there is exactly one arrow $\ii\to \ii'$.  

A functor 
$F:{\mathcal R}(w_0)\longrightarrow W-Cryst$ such that $F(\ii)={\mathcal
B}_\ii$ was defined in  \cite{bz3}. Moreover it was shown in  \cite{bz3} that
$\tilde \varphi^{\ii}_i=\varphi^{\ii'}_i\circ F\left(\ii\to
\ii'\right)$ for any $\ii ,\ii' \in  R(w_0)$. Therefore, the  
free $W$-crystal 
${\mathcal B}_{w_0}$ is well-defined and  is  equipped with the  
family of functions $\tilde \varphi_i:B_{w_0}\to \ZZ$, $i\in I$, such
that $\tilde \varphi_i(\tilde e_i(b))
=\tilde \varphi_i(b)+1$. 

\medskip

\noindent {\bf Remark}. Let $V_{\lambda^\vee}$ be the simple 
finite-dimensional 
$U_q({}^L\gg)$-module with the lowest weight $\lambda$ (e.g., 
$\lambda^\vee$ is an anti-dominant co-weight), and let 
${\mathcal B}(V_{\lambda^\vee})$ be the corresponding combinatorial 
crystal.
It follows from the results of Kashiwara that 
there is a morphism of $W$-crystals  
\begin{equation}
\label{eq:Kashiwara embedding}
f_{\lambda^\vee}:{\mathcal B}(V_{\lambda^\vee})\to \Lambda^\vee\times 
{\mathcal B}_{w_0}
\end{equation}
such that $\tilde \varphi_i(f_{\lambda^\vee}(b))=\sup\{n:b\in dom(\tilde e_i^{\,-n})\}$ for 
$b\in B(V_{\lambda^\vee})$, $i\in I$.

\end{document}